\documentclass[a4paper, 12pt]{article}

\usepackage{amsmath}
\usepackage{amsthm}
\usepackage{amsfonts}
\usepackage{amssymb}
\usepackage{latexsym}
\usepackage{color}
\usepackage{mathrsfs}
\usepackage{enumerate}
\usepackage{ulem, bm}

\usepackage{hyperref}
\hypersetup{colorlinks=true,urlcolor=black,linkcolor=black,citecolor=black}

\usepackage{mathtools}
\mathtoolsset{showonlyrefs=true}

\allowdisplaybreaks

\allowdisplaybreaks[4]

\usepackage{version}

\numberwithin{equation}{section}

\newtheorem{thm}{Theorem}[section]
\newtheorem{prop}[thm]{Proposition}
\newtheorem{lem}[thm]{Lemma}
\newtheorem{cor}[thm]{Corollary}
\newtheorem{rem}[thm]{Remark}
\newtheorem{definition}[thm]{Definition}

\newtheorem{exam}{Example}[section]

\newcommand{\form}{{\cal E}}
\newcommand{\dom}{{\cal F}}
\newcommand{\real}{{\mathbb R}}
\newcommand{\vareps}{\varepsilon}
\newcommand{\dis}{\displaystyle}

\def\sfA{{\sf A}}

\def\bE{{\mathbb{E}}}

\def\bM{{\mathbb{M}}}
\def\bN{{\mathbb{N}}}

\def\bP{{\mathbb{P}}}

\def\bR{{\mathbb{R}}}

\def\bfA{{\mathbf{A}}}

\def\cB{{\mathcal{B}}}

\def\cE{{\mathcal{E}}}
\def\cF{{\mathcal{F}}}
\def\cG{{\mathcal{G}}}

\def\cS{{\mathcal{S}}}

\newcommand{\ds}{\displaystyle}
\newcommand{\Capa}{{\rm Cap\,}}

\newcommand{\nest}{\mbox{\small {\sf N} \hspace*{-13.5pt} {\sf N}}}

\title{\sf
{Smooth Measures and Positive Continuous Additive Functionals Attached to a Compact Nest}}

\author{\sf
Takumu Ooi\thanks{Department of Mathematics,
Faculty of Science and Technology, Tokyo University of Science,
Noda, Chiba, 278-8510, Japan
({\sf ooitaku@rs.tus.ac.jp})},  \  Kaneharu Tsuchida\thanks{Department of Mathematics,
National Defense Academy,
Yokosuka, Kanagawa, 239-8686, Japan
({\sf tsuchida@nda.ac.jp}) } \ 
{\rm and}  \ Toshihiro Uemura\thanks{Department of Mathematics, 
Faculty of Engineering Science, 
Kansai University, Suita, Osaka, 564-8680, Japan
({\sf t-uemura@kansai-u.ac.jp})
}}

\date{}
\begin{document}

\maketitle 
\baselineskip=16pt 

\vspace*{-1cm}

\begin{abstract}
The relationship between smooth measures and positive continuous additive functionals is well known, and this correspondence is called the Revuz correspondence. We investigate the relationships between several types of convergence of smooth measures and convergence of positive continuous additive functionals, with a particular focus on the role of nests.  We provide conditions under which convergence of additive functionals implies convergence of the corresponding smooth measures. Our results cover convergence of smooth measures that are not  necessarily Radon, including nowhere Radon measures.
\end{abstract}

\section{Introduction}
In the theory of Markov processes, positive continuous additive functionals (PCAFs) and smooth measures play a crucial role. They are indispensable tools for transformations of Markov processes, such as killing and time change, and for analyzing their finer properties. The Revuz correspondence establishes a one-to-one relationship between PCAFs and smooth measures,  
serving as a crucial link between the analytic and probabilistic approaches.

The close connection between PCAFs and smooth measures is particularly evident in the analysis of transformed processes. For instance, the convergence of time-changed processes is closely related to the convergence of their associated PCAFs and smooth measures, a topic that has received significant attention in recent years (see, e.g., \cite{GRV16}, \cite{AK16}, \cite{Oo25}, and \cite{BK25+}).

Although the intrinsic and topological aspects of the well-known Revuz correspondence have been studied in specific cases, they have not been sufficiently explored within a general theoretical framework until recently. Intuitively, a smooth measure captures the sets that the process can hit because it does not assign a value to any set of zero capacity --- which, by definition, the process cannot hit.
A major challenge in defining a topological structure on the class of smooth measures is that a single smooth measure can behave differently depending on the choice of a ``nest''.  This class is also very broad, containing many non-Radon Borel measures.

In this paper, we investigate the connections between convergence of these two classes of objects, specifically focusing on providing conditions under which the convergence of a sequence of PCAFs is equivalent to, or implies, the convergence of their associated smooth measures under various topologies.

To establish our framework, we first outline the foundational settings and define the core objects of study. 
Let $E$ be a locally compact separable metric space and $m$ a positive Radon measure on $E$ with full topological support. 
We are given an $m$-symmetric Hunt process $\bM = (X_t, \mathbb{P}_x)$ on $E$. This is a strong Markov process whose sample paths are right-continuous with left limits, and quasi-left continuous, and its transition function is $m$-symmetric. 
We also assume that its associated Dirichlet form $(\cE, \cF)$ is regular. We denote by ${\cS}$ the class of smooth measures associated with $\bM$, and by ${\bfA_c^+}$ the class of PCAFs generated by $\bM$. Prior work has often focused on the more well-behaved subclass, ${\cS_0}$, the set of all smooth measures of finite energy integrals, to address the challenges posed by the broad nature of $\cS$.

Nishimori, Tomisaki, and the second and third authors, in \cite{NTTU25}, introduced a topology on the class ${\cal S}_0$. 
They defined a complete separable metric on this class that captures the energy of the measures, which allows these measures to be viewed as elements of the Dirichlet space.  
In this work, they also regarded the Revuz correspondence as a map (the Revuz map) from ${\cal S}$ to ${\bf A}^+_c$, and under certain assumptions, they showed its sequential compactness of the Revuz map. 
Later, in \cite{Oo26}, stronger results  were obtained, including: the continuity of the Revuz map from the class of measures restricted  to ${\cal S}_0$ to ${\bf A}^+_c$ with $L^p({\mathbb P}_x)$-norm for $1\le p<2$ by using the Fukushima decomposition of potentials of smooth measures, and the establishment of a homeomorphism of the Revuz map between ${\cal S}_0$ and the corresponding class of ${\bf A}^+_c$. 

These previous results are useful because proving the convergence of measures is often easier than proving the convergence of PCAFs. However, since they were mainly restricted to ${\cS_0}$, a broader framework is needed for the entire classes of ${\cS}$ and ${\bfA_c^+}$. To address the complexity of ${\cS}$ (which includes many non-Radon Borel measures) and provide a more stable framework, we introduce the notion of smooth measures and PCAFs attached to a nest $\{F_\ell \}$, which is an increasing sequence of compact sets such that ${\rm Cap}(K \setminus F_\ell)$,  the capacity of $K \setminus F_\ell$, goes to $0$ as $\ell \to \infty$ for any compact set $K \subset E$. 
For any smooth measure $\mu$,  there is a suitable nest $\{F_\ell\}$ such that its restriction $1_{F_\ell} \mu$ belongs to ${\cal S}_0$ for every $\ell$.  Motivated by these observations, we introduce the class ${\cal S}(\{F_\ell\})$ of smooth measures attached to a nest $\{F_\ell\}$. This class serves as a stepping stone because the choice of the nest plays a crucial role (see Example 3.1 below), thereby providing a more stable framework that allows us to work with measures whose supports are restricted to a nest,  as ${\cal S}$ itself  is equal to the union of the classes ${\mathcal S}(\{F_{\ell}\})$ over the family of nests. 

This paper explores several types of convergence on ${\cal S}(\{F_\ell\})$ for a fixed nest $\{F_\ell\}$, including weak and vague convergence on the nest,  to establish a suitable topological structure.  Although proving convergence of measures is often simpler, our primary motivation is the search for an appropriate concept of convergence in ${\cal S}(\{F_\ell\})$.  Therefore,  we primarily focus on conditions under which the convergence of PCAFs implies the convergence of measures.

Furthermore, we extend the topological approach of \cite{Oo26} (which established a homeomorphism of the Revuz correspondence on $\mathcal{S}_0$ in the $L^2$-sense) to a broader class of PCAFs within a unified $L^1$-framework. 
By incorporating the underlying reference measure $m$, the killing measure $\kappa$ from the Beurling--Deny decomposition, and the extended energy functional $\nu_0$ that characterizes the boundary ``absorption weight'' of the process, we propose a new, symmetric characterization of the Revuz correspondence in the $L^1$-sense. This formulation allows us to investigate a natural topology on the class of PCAFs with respect to the $L^1({\mathbb P}_{m+\kappa+\nu_0})$-norm. 
To keep the main text sharply focused on our central results concerning convergence on compact nests, we have placed this independent topological characterization in Appendix A, thereby ensuring a clearer and more concise presentation of the manuscript.

The rest of this paper is organized as follows. Section 2 establishes fundamental notions and notations related to Markov processes and Dirichlet forms. We then turn to Section 3, where we define the class of smooth measures attached to a nest, introduce the concepts of weak and vague convergence on a nest, and explore their related properties.  Section 4 moves on to introduce the weak convergence on a nest in the resolvent sense, providing sufficient conditions for deriving the convergence of smooth measures from the convergence of PCAFs. In Section 5, we provide illustrative examples. Notably, these include cases that previous research could not address, such as the convergence of Borel measures to a nowhere Radon measure and the approximation of singular Radon measures without finite energy integrals. Finally, as mentioned above, in Appendix A, we characterize the Revuz correspondence through the killing procedure and present sufficient conditions for the vague convergence of Radon smooth measures.

We close this introduction with a few words on notation:  We shall use $C_0(E)$ to denote the set of continuous functions with compact support in $E$ and denote by $\cB(E)$
the set of all Borel measurable functions on $E$. If $\cG$ is a set of functions, then $\cG_b$ (resp. $\cG_+$) denotes the set of bounded (resp. non-negative) functions in $\cG$.

\section{Preliminaries}
Let $E$ be a locally compact separable metric space and $m$ a positive Radon measure on $E$
with full topological support. Define $E_{\partial} = E \cup \{\partial\}$ as the one-point compactification of $E$.

Let ${\mathbb M} = (\Omega,\cF,\cF_t,X_t,\bP_x,\zeta)$ be an $m$-symmetric Hunt process on $E$. 
Here $\{\cF_t\}$ is the minimal augmented filtration and $\zeta=\inf\{t\ge 0 : X_t=\partial\}$ is the lifetime 
of the process.

The transition semigroup of ${\mathbb M}$ is denoted by $\{p_t\}_{t \ge 0}$, defined as $p_t f(x) ={\mathbb E}_x[f(X_t)]$. 
According to \cite[Lemma 1.4.3]{FOT11}, this semigroup uniquely determines a strongly continuous 
Markovian semigroup $\{T_t\}_{t \ge 0}$ on $L^2(E; m)$. 
The Dirichlet form on $L^2(E;m)$ associated with $\bM$ is  then given by
\begin{align}
 \begin{cases}
  \ds \cE(u,v) = \lim_{t \downarrow 0} \frac{1}{t} (u - T_t u, v), \medskip\\
  \ds \cF = \left\{u \in L^2(E;m) : \lim_{t \downarrow 0} \frac{1}{t}
  (u - T_t u, u) < \infty\right\}, 
 \end{cases}
\end{align}
where $(f,g)$ and $\|f\|_2$ denote the $L^2$-inner product and $L^2$-norm, respectively.

For $\alpha > 0$, we define the form $\form_\alpha$ by 
\begin{align}
 \cE_{\alpha}(u,v) := \cE(u,v) + \alpha (u,v). 
\end{align}
It is known that $\cF$ is a real Hilbert space with
the inner product $\cE_1$. We denote the corresponding $\form_1$-norm by $\|\cdot\|_{\cE_1}$, defined as 
\begin{align}
 \|u\|_{\cE_1} & :=\sqrt{\cE_1(u)} := \sqrt{\cE_1(u,u)}, \quad u \in \cF.
\end{align}

Throughout this paper, we always assume that $(\cE,\cF)$ is regular. This means that 
 $\cF \cap C_0(E)$ is dense in $\cF$ with respect to the $\cE_1$-norm and  also dense in $C_0(E)$ 
 with respect to the sup-norm $\|\cdot\|_{\infty}$.

 Next we define the capacity associated with $(\cE,\cF)$ by \(\Capa(O) :=\ds \inf\big\{\cE_1(u,u) : u \in \cF : u \ge 1\ m\text{-a.e. on } O\big\}\) for an open set $O \subset E$, where the infimum of the empty set is \(\infty\), and \(\Capa(A) := \inf\left\{\Capa(O) : A \subset O,\ O\ \text{is open}\right\}\) for any set $A \subset E$.

We call a set $A \subset E$ {\it exceptional} if $\Capa(A) = 0$.
A statement depending on $x \in A$ is said to hold  {\it quasi-everywhere} ({\it q.e.})  on $A$
if there exists an exceptional set $N \subset A$ such that the statement holds for every 
$x \in A \setminus N$. A \([-\infty, \infty]\)-valued function $u$ defined q.e. on $E$ is called {\it quasi-continuous} if, 
for any $\varepsilon > 0$, there exists an open set $G \subset E$ such that $\Capa(G) < \varepsilon$ 
and the restriction $u|_{E \setminus G}$ of $u$ to $E\setminus G$ is finite and continuous on $E\setminus G$.  
By \cite[Theorem 2.1.3]{FOT11}, every function $u \in \cF$ admits a quasi-continuous version. 
Therefore, in the rest of this paper, we always assume that any function $u \in \cF$ is quasi-continuous.

We call an increasing sequence of 
closed sets $\{F_\ell\}$ a {\it generalized nest} if 
\begin{equation}  
{\rm Cap}(K\setminus F_\ell) \to 0 \quad {\rm as} \ \ell \to \infty \quad {\rm for \ any \  compact  \ set} \label{nest}
\  K\subset E. 
\end{equation}
A positive Borel measure $\mu$ on $E$ is said to be
{\it smooth} if it satisfies the following conditions:
\begin{itemize}
 \item $\mu$ charges no set of zero capacity.
 \item There exists a 
       generalized nest $\{F_{\ell}\}$ such that
       $\mu(F_\ell)<\infty$ for each $\ell$.
       
\end{itemize}
Denote by $\cS$ the class of all smooth measures. 

A positive Radon measure $\mu$ on $E$ is called a {\it measure
of finite energy integral} if there exists a constant
{
$C > 0$ 
} 
such that 
\begin{equation}
\int_E |v(x)| \mu(dx) \le C \sqrt{\form_1(v,v)}, \quad {\rm for \ any \ }v\in \dom\cap C_0(E). \label{eq:def_S_0}
\end{equation}

Let $\cS_0$ be the family of all measures of finite energy integrals.
From  \cite[p.84]{FOT11}, we know that $\cS_0$ is a subclass of $\cS$. 
For any $\mu \in \cS_0$ and $\alpha > 0$, the Riesz representation theorem states that 
there exists a unique function $U_{\alpha} \mu \in \cF$ such that
\begin{align}
 \cE_{\alpha}(U_{\alpha}\mu,v) = \int_{E} v(x) \mu(dx),
 \quad \text{for any $v \in \cF \cap C_0(E)$}.\label{eq:def_U_1}
\end{align}
We call $U_{\alpha} \mu$ an {\it $\alpha$-potential} of $\mu$.  When $(\form, \dom)$ is transient, we also define \(\mathcal{S}_0^{(0)}\) as the set of positive Radon measures satisfying \eqref{eq:def_S_0} with \(\mathcal{E}_1\) replaced by 
\(\mathcal{E}\), and, for \(\mu \in \mathcal{S}_0^{(0)}\), there exists \(U\mu \in \mathcal{F}_e\) such that \eqref{eq:def_U_1} with \(\mathcal{E}_{\alpha}\) replaced by \(\mathcal{E}\). According to \cite[Theorem 2.2.4]{FOT11},  for any smooth measure $\mu$, there exists a generalized 
nest $\{F_\ell\}$ such that $1_{F_\ell} \mu \in {\cal S}_0$ for each $\ell$. Furthermore, the generalized nest \(\{F_\ell\}\) satisfying this condition can be chosen as an increasing sequence of compact sets, which is called a {\it generalized compact nest}. For this reason, throughout this paper, we will simply call such a generalized compact nest a {\it nest}.

Next we introduce the definition of a positive continuous additive functional.
\begin{definition}
 \label{def-PCAF}
A function \(\sfA : \Omega \times [0,\infty) \to [0,\infty]\) is a positive continuous additive functional (PCAF) if the following conditions hold:
 \begin{enumerate}
 \item for each $t \ge 0$, $\sfA_t(\cdot)$ is $\cF_t$-measurable,
  \item
       there exist a set $\Lambda \in \cF_{\infty}$ and
       an exceptional set $N \subset E$ such that $\bP_x(\Lambda) = 1$
       for any $x \in E \setminus N$ and $\theta_t \Lambda \subset \Lambda$
       for any $t > 0$, and moreover,
       for each $\omega \in \Lambda$, ${\sf A}_0(\omega) = 0$,
       $|\sfA_t(\omega)| < \infty$ for any $t < \zeta(\omega)$,
       $\sfA_t(\omega) = {\sf A}_{\zeta(\omega)}(\omega)$ for any
       $t \ge \zeta(\omega)$ and
       \begin{align}
	{\sf A}_{t + s}(\omega) = {\sf A}_s(\omega) + {\sf A}_t(\theta_s \omega),\
	\text{for any $t,s \ge 0$},
       \end{align}
  \item
       for each $\omega \in \Lambda$, the map
       $t \mapsto {\sf A}_t(\omega)$ is continuous, where $\Lambda$ is defined in the condition 2.
 \end{enumerate}
\end{definition}
The set of all PCAFs is denoted by $\bfA_c^+$. Two PCAFs $\sfA^{(1)}$ and $\sfA^{(2)}$, 
 are said to be {\it equivalent} if for each $t > 0$, $\bP_x(\sfA_t^{(1)} = \sfA_t^{(2)}) = 1$ q.e. $x \in E$.
If $\sfA^{(1)}$ and $\sfA^{(2)}$ are equivalent, we write $\sfA^{(1)} \sim \sfA^{(2)}$. 
This relation $\sim$ clearly  defines an equivalence relation on $\bfA_c^+$.
The following theorem describes  the concrete relationship  between $\cS$ and
$\bfA_c^+$.

\begin{thm}[{see \cite[Theorem 5.1.3, 5.1.4]{FOT11}}]
 \label{thm-Revuz}
 The quotient space of $\bfA_c^+$ under the equivalence relation
 $\sim$ and the family $\cS$ are in one-to-one correspondence under
 the following relation:
 For $\mu \in \cS$ and $\sfA \in \bfA_c^+$,
 \begin{align} \label{eq-Revuz}
  \mathbb{E}_{hm} \left[ \int_0^\infty e^{-\alpha s} f(X_s) \, d\sfA_s \right] 
  = \int_E f(x) R_\alpha h(x)\, \mu(dx)
 \end{align}
for any $f, h \in \cB_+(E)$. Moreover, $(\ref{eq-Revuz})$ is equivalent to 
\begin{equation}
 \mathbb{E}_{hm} \left[ \int_0^t  f(X_s) \, d\sfA_s \right] =
  {
   \int_{E} \left(\int_{0}^{t} p_s h(x)\, ds\right) f(x)\, \mu(dx)}
\label{eq-Revuz2}
\end{equation}
 for any $f, h \in \cB_+(E)$ and \(t\geq 0\).
\end{thm}
Hence these equations \eqref{eq-Revuz} and \eqref{eq-Revuz2} are called the {\it Revuz correspondence}. 

\section{Smooth Measures Attached to a Compact Nest and Weak Convergence}

The class ${\cal S}$ is quite large (it is larger than the set of all Radon measures charging 
no sets of zero capacity, and includes many singular Borel measures),  and so is ${\bf A}_c^+$. 

In order to manage this larger class, we introduce a class of smooth measures attached to a nest 
in this paper.

Let $\{F_\ell\}$ be a nest. A smooth measure $\mu$ is said to be \textit{attached 
to  $\{F_\ell\}$} if the restriction $1_{F_\ell} \mu$ belongs to ${\cal S}_0$ for each $\ell$. Denote by 
${\cal S}(\{F_\ell\})$ the set of all smooth measures attached to the nest $\{F_\ell\}$. 
Since the set $F_\ell$ is compact for each $\ell$ and 
$(\form, \dom)$ is regular, there exists $\varphi \in \dom \cap C_0(E)$ such that $\varphi \ge 0$ 
and $\varphi=1$ on $F_\ell$. Then, for $\mu \in {\cal S}(\{F_\ell\})$, it follows from the definition 
of ${\cal S}_0$ that 
\begin{equation} \label{eq:Radon}
\mu(F_\ell)  \le \int_{F_\ell} \varphi (x) \mu(dx) =\int_E \varphi(x) (1_{F_\ell}\mu)(dx) \le 
C_\ell \sqrt{\form_1(\varphi, \varphi)}<\infty
\end{equation}
holds for each $\ell$. It follows from the fact that \(E\) is locally compact and separable that
\begin{equation}
   \label{eq:smooth_nest}
\mu \Big(E\setminus \bigcup_{\ell=1}^\infty F_\ell\Big)=0.
\end{equation}
We also see that 
$$
{\cal S}_0 \subset {\cal S}(\{F_\ell\}).
$$ 
We call a nest $\{F_\ell\}$ an {\it exhaustion} of $E$ if $\bigcup_{\ell=1}^\infty F_\ell =E$ and, for any compact set \(K\), there exists \(N\in \mathbb{N}\) such that \(K \subset F_N.\) Then \eqref{eq:Radon} implies that a measure $\mu \in{\cal S}(\{F_\ell\})$ is Radon.

Given another nest $\{G_\ell\}$,  it follows that $\{F_\ell \cap G_\ell\}$ is also a nest.  In fact,  for any compact set $K\subset E$, 
\begin{align*}
{\rm Cap}(K\setminus (F_\ell \cap G_\ell)) 
& ={\rm Cap}\Big( \big(K\setminus F_\ell \big) \cup 
\big(K\setminus G_\ell\big)\Big) \\
& \le {\rm Cap}(K\setminus F_\ell) + 
{\rm Cap}(K\setminus G_\ell) \to 0 \ \ {\rm as} \ \  \ell \to\infty.
\end{align*}
Moreover both $\mu \in {\cal S}(\{F_\ell\})$ and $\nu \in {\cal S}(\{G_\ell\})$  belong to ${\cal S}(\{F_\ell \cap G_\ell\})$. 
 Hence, ${\cal S}(\{F_\ell\}) \cup {\cal S}(\{G_\ell\}) \subset {\cal S}(\{F_\ell \cap G_\ell\})$.  However, the reverse inclusion does not generally hold (see Example \ref{ex-01} below).

Let \nest \ be the set of all nests.  Then by virtue of \cite[Theorem 2.2.4]{FOT11},  we see that 
$$
{\cal S}= \bigcup_{\{F_\ell\} \, \in \, {\nest} } {\cal S} (\{F_\ell\}).
$$

In the following we fix a nest $\{F_\ell\} \in \nest$\,
and define a function $\rho$ on ${\cal S}(\{F_\ell\}) \times \mathcal{S}(\{F_{\ell}\})$: 
\begin{equation} \label{metric}
\rho(\mu, \nu):= \sum_{\ell=1}^\infty  \frac 1{2^{\ell}} \Big( 1\wedge 
\sqrt{\form_1\big(U_1(1_{F_\ell} \mu)-U_1(1_{F_\ell} \nu)\big)} \, \Big).
\end{equation}
We have defined a complete and separable metric $\rho_0$ on ${\cal S}_0$ by using the 
$1$-potentials  in \cite{NTTU25}:
$$
\rho_0(\mu, \nu):=\sqrt{\form_1(U_1\mu-U_1\nu)}, \quad \mu, 
\nu \in {\cal S}_0.
$$
Then $\rho$ is also written as follows:
\begin{equation} \label{rho}
\rho(\mu, \nu) =  \sum_{\ell=1}^\infty  \frac 1{2^\ell} 
\Big( 1\wedge \rho_0(1_{F_\ell}\mu, 1_{F_\ell}\nu)\Big).
\end{equation}
Note that although $\rho$ defines a metric on $\mathcal{S}(\{F_\ell\})$
as we will see in Lemma \ref{lem:metric},
the metric space $(\mathcal{S}(\{F_\ell\}), \rho)$
is not complete in general, unlike the metric $\rho_0$ on $\mathcal{S}_0$.
We will give a specific counterexample and
discuss the reason for this incompleteness in Remark \ref{rem:sep} below.

\begin{lem}
  \label{lem:metric}
  The function $\rho$ defines a separable metric on ${\cal S}(\{F_\ell\})$. 
\end{lem}
\begin{proof}
  Let $\mu, \nu, \xi \in \mathcal{S}(\{F_{\ell}\})$.
  The symmetric property $\rho(\mu,\nu) = \rho(\nu,\mu)$ and the triangle inequality
  $\rho(\mu, \nu) \le \rho(\mu,\xi) + \rho(\xi,\nu)$ follows immediately from the fact
  that $\rho_{0}$ is a metric on $\mathcal{S}_{0}$.

  Now we show that $\rho(\mu, \nu) = 0$ implies $\mu = \nu$.
  Assume that $\rho(\mu,\nu) = 0$. By the definition of $\rho$,
  we have $\rho_{0}(1_{F_{\ell}}\mu, 1_{F_{\ell}}\nu) = 0$ for all $\ell \ge 1$.
  Since $\rho_{0}$ is a metric on $\mathcal{S}_{0}$,
  this means $1_{F_{\ell}} \mu = 1_{F_{\ell}}\nu$ as
  measures in $\mathcal{S}_{0}$ for each $\ell \ge 1$. 
  By \eqref{eq:smooth_nest}, any measure in $\mathcal{S}(\{F_{\ell}\})$
  carries no mass on $E \setminus \bigcup_{\ell=1}^{\infty }F_{\ell}$.
  Therefore, both $\mu$ and $\nu$ vanish on
  $E \setminus \bigcup_{\ell=1}^{ \infty}F_{\ell}$.
  Since $\mu$ and $\nu$ coincide on each $F_{\ell}$ and have no mass
  outside the union of the nest, we conclude that $\mu = \nu$ on $E$.

We define a map \(\iota: \mathcal{S}(\{F_\ell\}) \to \mathcal{S}_0^{\mathbb{N}}\) by \(\iota(\mu):=(1_{F_\ell}\mu)_\ell\). The argument above showing that \(\rho(\mu,\nu)=0\) implies \(\mu=\nu\) also shows that \(\iota\) is injective.
For \((\mu_k)_k, (\nu_k)_k \in \mathcal{S}_0^\mathbb{N}\), we set 
\[\rho_0^\mathbb{N}\bigl((\mu_k)_k, (\nu_k)_k\bigr) :=
\sum_{k=1}^{\infty}\frac{1}{2^k}
\Big( 1\wedge \rho_0(\mu_k, \nu_k)\Big).\]
Then \((\mathcal{S}(\{F_\ell\}),\rho)\) is isometric to \((\iota(\mathcal{S}(\{F_\ell\})), \rho_0^\mathbb{N})\). Since \((\mathcal{S}_0,\rho_0)\) is separable (\cite[Lemma 3.4]{NTTU25}), the countable product \((\mathcal{S}_0^{\mathbb{N}}, \rho_0^\mathbb{N})\) is also separable, and, since every subspace of a separable metric space is separable, \((\mathcal{S}(\{F_\ell\}),\rho)\) is separable.
\end{proof}

\begin{rem}
 \label{rem:sep}
 \begin{enumerate}[(1)]
 \item It should be noted that the metric space $(\mathcal{S}(\{F_{\ell}\}), \rho)$
   is not complete in general. For instance, consider the one-dimensional
   Brownian motion. If we set $F_{\ell} = [-\ell,\ell]$ and $\mu_{n} = \delta_{1+1/n}$,
   the sequence $\{\mu_{n}\}$ forms a $\rho$-Cauchy sequence.
   Indeed, for $\ell = 1$, we have $1_{F_{1}} \mu_{n} = 0$ for all $n \ge 1$.
   For $\ell \ge 2$, the point $1 + 1/n$ is eventually inside $F_{\ell}$, so
   $1_{F_{\ell}} \mu_{n} = \delta_{1 + 1/n}$ and \(\delta_{1 + 1/n}\) converges to \(\delta_1\) in $\rho_{0}$ as $n \to \infty$ by the continuity of the resolvent density. However, if it were to converge to some measure $\mu \in \mathcal{S}(\{F_\ell\})$, 
   the convergence would require $1_{F_1}\mu = 0$
   and $1_{F_\ell}\mu = \delta_1$ for $\ell \ge 2$. 
   This implies $\mu = \delta_1$ on the whole space $E$, which forces 
   $1_{F_1}\mu = \delta_1 \neq 0$, leading to a contradiction. 
   Thus, it does not converge to any measure in $\mathcal{S}(\{F_{\ell}\})$.
 \item The reason for this incompleteness is
   that the discontinuous indicator function $1_{F_\ell}$ cannot properly capture
   the mass moving across the boundary $\partial F_\ell$.
   If one wishes to make the space complete,
   one could assume that the nest $\{F_\ell\}$ satisfies
   $F_\ell \subset F_{\ell+1}^o$ and, by the regularity of the Dirichlet form, we take a sequence of
   continuous cut-off functions $\{\psi_\ell\} \subset \mathcal{F} \cap C_0(E)$
   such that $0 \le \psi_\ell \le 1$, $\psi_\ell = 1$ on $F_\ell$,
   and $\psi_\ell = 0$ on $F_{\ell + 1}^{c}$.
   By replacing $1_{F_\ell}$ with $\psi_\ell$, the metric $\rho$ can be redefined as
   \begin{align*}
     \widetilde{\rho}(\mu,\nu)
     := \sum_{\ell=1}^{\infty } \frac{1}{2^{\ell}}
     \left( 1 \wedge \rho_{0}(\psi_{\ell}\mu, \psi_{\ell}\nu)\right).
   \end{align*}
   Indeed, for a \(\widetilde{\rho}\)-Cauchy sequence \(\{\mu_n\}_n\), there exists \(\{\mu^{(\ell)}\}_\ell \subset \mathcal{S}_0\) such that \(\rho(\psi_\ell \mu_n, \mu^{(\ell)}) \to 0\) as \(n\to \infty\). Since \(\psi_\ell f \in \mathcal{F}\) for any \(f\in \mathcal{F}\cap C_0(E)\), we have 
\begin{align*}
  \lefteqn{\int_E f \psi_\ell \,d\mu^{(\ell +1)} = \mathcal{E}_1(U_1\mu^{(\ell+1)}, f\psi_\ell) =  \lim_{n\to \infty} \mathcal{E}_1(U_1(\psi_{\ell+1} \mu_n), f\psi_\ell)}\\
   &= \lim_{n\to \infty} \int_E f \psi_{\ell +1} \psi_\ell \,d\mu_n = \lim_{n\to \infty} \int_E f  \psi_\ell \,d\mu_n = \lim_{n\to \infty} \mathcal{E}_1(U_1(\psi_{\ell} \mu_n), f) = \int_E f \,d\mu^{(\ell)}
\end{align*}   
for \(f\in \mathcal{F}\cap C_0(E)\). By the regularity, this equality holds for any \(f \in C_0(E)\) and so \(\psi_{\ell}\mu^{(\ell+1)} = \mu^{(\ell)}\) holds. We define \(\mu(B):=\lim_{\ell \to \infty}\mu^{(\ell)}(B\cap F_\ell)\) for a Borel set \(B\), it holds that \(\mu \in \mathcal{S}(\{F_\ell\})\) and \(\mu_n \to \mu\) in \(\widetilde{\rho}\).

   In this paper, however, we proceed with the current definition
   using $1_{F_\ell}$, as it is the most natural and powerful tool
   for our main purpose: to evaluate the convergence
   of singular measures attached to suitably chosen nests,
   which will be demonstrated in Section \ref{sec_example}.
   Note that an appropriate subsequence of an exhaustion also satisfies
   the condition $F_\ell \subset F_{\ell+1}^o$.
 \end{enumerate}
\end{rem}

\begin{exam}  \label{ex-01}
Let  $(\form,\dom)=(\frac 12 {\mathbb D}, H^1(\real^d)), \, d\ge 2$. 
Take $a\in\real^d\setminus\{0\}$ and $\beta\in \real$. Set  $\mu(dx):=|x|^{-\beta}dx$, $\nu(dx):=|x-a|^{-\beta}dx$  and \(\bar{B}(y;r):=\{x\in \mathbb{R}^d : |x-y| \leq r\}\) for \(y\in \mathbb{R}^d\) and \(r>0\). Then by \cite[Example 5.1.1]{FOT11}, it is known that $\mu$ and 
$\nu$ are smooth measures.  For $\ell \in{\mathbb N}$,  set 
$$
F_\ell :=\{ x \in \mathbb{R}^d : 1/\ell \le  |x| \le  \ell  \} \qquad  {\rm and} \qquad G_\ell :=\{ x \in \mathbb{R}^d :  1/\ell \le  |x-a| \le \ell  \}.
$$
Assume that $\beta \ge  (d+2)/2$ holds. Then $\mu\in{\cal S}(\{F_\ell\})$ and 
$\nu \in{\cal S}(\{G_\ell\})$ but $\mu \not\in {\cal S}(\{G_\ell\})$ and 
 $\nu \not\in {\cal S}(\{F_\ell\})$. Indeed, for \(\ell >2|a| \vee 2|a|^{-1}\), since \(G_\ell \supset \bar{B}(0;|a|/2)\) and 
 \( \bar{B}(0;|a|/2)\times \bar{B}(0;|a|/2) \supset  \{(x,y) \in \mathbb{R}^d \times \mathbb{R}^d : |x|\leq |a|/6, 2|x| \leq |y| \leq 3|x| \}=:K \), for \(d\geq 3\), we have
\begin{eqnarray*}
\iint_{G_\ell\times G_\ell} g_1(x,y) |x|^{-\beta} |y|^{-\beta} dxdy & \geq &C\iint_{K} |x-y|^{2-d} |x|^{-\beta} |y|^{-\beta} dxdy\\
\geq C\int_{\{|x|\leq \frac{|a|}{6}\}} |x|^{2-d} |x|^{-\beta} \int_{\{2|x|\leq |y| \leq {3|x|}\}} |y|^{-\beta} dy dx & \geq &C \int_0^{\frac{|a|}{6}} r^{d+1-2\beta} dr = \infty,
\end{eqnarray*} 
and so \(1_{G_\ell}\mu \not \in \mathcal{S}_0\) by \cite[Exercise 4.2.2]{FOT11}. Here \(g_1\) is a \(1\)-resolvent density for Brownian motion. We remark that, for \(\beta <(d+2)/2\) and \(d\geq 3\), by letting \(B:=\bar{B}(0;R)\times \bar{B}(0;R)\), we have
\begin{eqnarray*}
\lefteqn{ \iint_{B}\frac{g_1(x,y)}{|x|^{\beta} |y|^{\beta}} dxdy}\\
&=& 2\iint_{\{(x,y)\in B : |y|< \frac{|x|}{2}\}} \frac{g_1(x,y)}{|x|^{\beta} |y|^{\beta}}  dxdy +\iint_{\{(x,y)\in B : \frac{|x|}{2} \leq |y| \leq 2|x|\}} \frac{g_1(x,y)}{|x|^{\beta} |y|^{\beta}} dxdy \\
&\leq & C \int_{\bar{B}(0;R)} |x|^{2-d-\beta} \int_{\bar{B}(0;\frac{|x|}{2})} |y|^{-\beta}dy dx +\iint_{\{(x,y)\in B : |x-y| \leq 3|x|\}} |x-y|^{2-d} |x|^{-2\beta}dydx\\
&=& C \int_{\bar{B}(0;R)} |x|^{2-2\beta} dx + C \int_{\bar{B}(0;R)} |x|^{2-2\beta} dx \ <\ \infty,
\end{eqnarray*}
and so \(\mu \in \mathcal{S}(\{G_\ell\}).\) The case of \(d=2\) is similar. We also see $\mu, \nu \in {\cal S}(\{F_\ell \cap G_\ell\})$.  
 Furthermore, the smooth measure $\mu + \nu$ belongs to
 $\cS(\{F_{\ell} \cap G_{\ell}\})$, but it belongs to neither $\cS(\{F_{\ell}\})$
 nor $\cS(\{G_{\ell}\})$. Therefore $\cS(\{F_{\ell} \cap G_{\ell}\}) \not\subset
 \cS(\{F_{\ell}\}) \cup \cS(\{G_{\ell}\})$. This is a reason why we consider convergence of smooth measures attached to a fixed nest. In particular, define 
$\mu_n(dx):=|x|^{-\frac d2 -1+\frac 1n} dx \in{\cal S}(\{F_\ell\}) \cap {\cal S}(\{G_\ell\})$ for any 
$n\in{\mathbb N}$.  In this case,  we find that $\{\mu_n\}$ is $\rho$-Cauchy in ${\cal S}(\{F_\ell\})$ by the monotonicity and the limit $\mu(dx):=|x|^{-\frac d2 -1} dx \in{\cal S}(\{F_\ell\})$
 but not in ${\cal S}(\{G_\ell\})$ since $\mu \not\in {\cal S}(\{G_\ell\})$.
\end{exam}

 Smooth measures are not necessarily Radon measures or finite measures, and so we cannot consider the vague convergence or the weak convergence for smooth measures. We now introduce a notion of convergence of smooth measures, which we call {\it weak convergence on a nest}.

\begin{definition}
 Let $\{F_\ell\}$ be a nest.  A sequence of smooth measures $\{\mu_n\}$ is said to 
\textit{converge to a measure $\mu$ weakly on the nest $\{F_\ell\}$} if 
$\{\mu_n\} \subset {\cal S}(\{F_\ell\})$ and  
\begin{equation} \label{weak1}
\lim_{n\to \infty} \int_{F_\ell} f(x) \mu_n(dx)=\int_{F_\ell} f(x) \mu(dx), \quad {\rm for \ every} \ 
f\in C(\{F_\ell\}) \ {\rm and} \ \ell \in{\mathbb N}.
\end{equation}
Here $C(\{F_\ell\})$ denotes the set of all extended real-valued functions $f$ defined q.e. on $E$ 
such  that the restriction of $f$ to $F_\ell$ is continuous for each $\ell$ (see \cite[\S 2.1]{FOT11}).
\end{definition} 
 
It is crucial to note that, even if the sequence $\{\mu_n\}$ converges weakly to $\mu$ on the nest $\{F_\ell\}$, the limit measure $\mu$ may not belong to the class ${\cal S}(\{F_\ell\})$.  For instance, 
in Example \ref{ex-01}, the sequence of measures $\nu_n(dx):=(|x-a|^{-\beta}\wedge n )dx$ 
converges weakly to the smooth measure $\nu(dx)=|x-a|^{-\beta}dx$ on the nest $\{ F_\ell\}$ for any $\beta<d$. 
Despite this convergence, $\nu$ is not a member of ${\cal S}(\{F_\ell\})$ when  $ (d+2)/2\le \beta<d$. 
On the other hand, for any $\beta$, the sequence $\{\nu_n\}$ converges weakly to $\nu$ on the nest 
$\{G_\ell\}$, and in this case, both $\{\nu_n\}$ and $\nu$ are members of  ${\cal S}(\{G_\ell\})$. In general, it does not automatically follow that the limiting measure is also smooth.  For example, in the case of Brownian motion on \(\mathbb{R}^2\), we set \(d\mu_n := n^2 \varphi(nx)dx \) for \(\varphi \in C_0(\mathbb{R}^2)\) satisfying \(\int \varphi dx =1\) and \(0\leq \varphi \leq 1\), and \(F_\ell := \{x\in \mathbb{R}^2 :|x|\leq \ell\} \). Then \(\mu_n \in \mathcal{S}(\{F_\ell\})\) and \(\mu_n\) converges weakly to the Dirac measure \(\delta_0\) at \(0\) on the nest \(\{F_\ell\}\), but \(\delta_0\) is not smooth. See \cite{F99, H04} for sufficient conditions for a measure to be smooth.

\begin{lem} Let $\{\mu_n\}\subset {\cal S}(\{F_\ell\})$. Assume that 
$\dis  M_\ell:=\sup_{n\in{\mathbb N}} \form_1(U_1(1_{F_\ell} \mu_n), U_1(1_{F_\ell} \mu_n))<\infty
$ for each $\ell\in{\mathbb N}$. Then 
$$
\sup_n \mu_n(F_\ell)<\infty
$$ 
for each $\ell$.  In particular, 
$$
\sup_n \mu_n(K)<\infty
$$ 
holds for any compact set $K$, provided that $\{F_\ell\}$ is an exhaustion of $E$. 
\end{lem}
 \begin{proof}
  Take a $\varphi \in \dom\cap C_0(E)$ satisfying $\varphi=1$ on $F_\ell$ and $0\le \varphi \le 1$ 
by means of the regularity. Then, for each $\ell$  we see from the Schwarz inequality that 
\begin{align*}
\mu_n(F_\ell)  & \le \int_{F_\ell} \varphi d\mu_n =\form_1(U_1(1_{F_\ell}\mu_n), \varphi) \\
& \le \sqrt{\form_1(U_1(1_{F_\ell}\mu_n),U_1(1_{F_\ell}\mu_n))}\sqrt{\form_1(\varphi, \varphi)}
\le \sqrt{M_\ell}\sqrt{\form_1(\varphi, \varphi)} <\infty. 
\end{align*}
The second assertion holds provided that $\{F_\ell\}$ is an exhaustion of $E$. In this case, for any compact set 
$K\subset E$, there exists an $\ell$ such that $K\subset F_\ell$, whence we see$$
\sup_n \mu_n(K) \le \sup_n \mu_n(F_\ell) <\infty. 
  $$
  {
  The proof is completed.}
 \end{proof}

\begin{prop} \label{prop:weak} 
  Let $\{F_\ell\}$ be a nest. Let $\{\mu_n\} \subset \mathcal{S}(\{F_\ell\})$ and $\mu \in \mathcal{S}(\{F_\ell\})$.  
   If  $\dis \lim_{n \to \infty} \rho(\mu_n, \mu) = 0$, then the sequence  $\{\mu_n\}$ converges to $\mu$ weakly on the nest $\{F_\ell\}$.
  Furthermore, if the nest $\{F_\ell\}$ is an exhaustion of $E$, then the sequence \(\{\mu_n\}\) converges vaguely to \(\mu\), that is,
  $$
  \lim_{n\to \infty} \int_{E} f(x) \mu_n(dx)=\int_{E} f(x) \mu(dx), \quad {\rm for \ every} \ f\in C_0(E).
  $$
\end{prop}

\begin{proof} 
 Let $ f \in C(\{F_\ell\})$. Take any $\ell \in{\mathbb N}$ and fix it. Since $F_\ell$ is a compact 
 set and $ f$ is continuous on $F_\ell$, we can find a function $g\in C_0(E)$ with 
 $f=g$ on $F_\ell$.  By the regularity, \eqref{eq:Radon} and the previous lemma, for each 
 $\vareps>0$, there exists $h \in \dom\cap C_0(E)$ such that 
 \[\| g-h\|_\infty < \frac{\vareps}{2} \big(\max\{\sup_n \mu_n(F_\ell), \mu(F_\ell)\} +1\big).\] Then 
\begin{align*}
& \!\!\!\!\! \Big| \int_{F_\ell}  f(x)\mu_n(dx) - \int_{F_\ell}   f(x)\mu(dx) \Big| 
 =  \Big| \int_{F_\ell} g(x)\mu_n(dx) - \int_{F_\ell}  g(x)\mu(dx) \Big| \\
& \le  \Big| \int_{F_\ell} g(x)\mu_n(dx) - \int_{F_\ell} h(x)\mu_n(dx) \Big| + 
\Big|\int_{F_\ell} h(x)\mu_n(dx)  - \int_{F_\ell} h(x)\mu(dx)  \Big|  \\
& \quad  + \Big| \int_{F_\ell} h(x)\mu(dx)-\int_{F_\ell} g(x)\mu(dx) \Big| \\
& \le \|g-h\|_{\infty} \mu_n(F_\ell) + \big| \form_1(U_1(1_{F_\ell}\mu_n), h) - 
\form_1(U_1(1_{F_\ell}\mu), h) \big|  + \|g-h\|_{\infty} \mu(F_\ell) \\
& < \vareps + \big| \form_1(U_1(1_{F_\ell}\mu_n), h) - 
\form_1(U_1(1_{F_\ell}\mu), h) \big| 
\end{align*}
holds. Thus tending $n\to\infty$ and then $\vareps \to 0$, the right hand side goes to $0$.

When $\{F_\ell\}$ is an exhaustion of $E$,  for any $ f \in C_0(E)$, there exists an $\ell$ such that ${\rm  supp}[ f] \subset F_\ell$.  Consequently, we have $\int_E  f d\mu_n = \int_{F_\ell}  f d\mu_n$ and $\int_E  f d\mu = \int_{F_\ell}  f d\mu$. The required convergence $\lim_{n \to \infty}\int_E f d\mu_n = \int_E  f d\mu$ then follows directly from the first part of the proof by choosing this $\ell$. 
 \end{proof}

\medskip

In the rest of this section, we summarize some results on algebraic properties for 
\(\mathcal{S}(\{F_{\ell}\})\) and the convergence of PCAFs based on their smooth measures in 
\(\mathcal{S}(\{F_{\ell}\})\), which were established for \(\mathcal{S}_0\) in \cite{NTTU25}.
For a PCAF \({\sf A }\) and a  Borel set $B$, we define \((1_B{\sf A})_t:=\int_0^t 1_B(X_s) d{\sf A}_s\).

\begin{prop}[{\it c.f.} {\cite[ Proposition 3.10]{NTTU25}}]\label{smoooth_alg} The following 
properties hold for ${\cal S}(\{F_\ell\})$. 
\begin{itemize}
\item[\((1)\)] {\sf (Monotonicity)} Let $\mu$ and $\nu$ be any Borel measures on $(E,{\cal B}(E))$. Assume that $\mu \le \nu$ and $\nu \in{\cal S}(\{F_\ell\})$.  Then $\mu \in{\cal S}(\{F_\ell\})$ and $U_\alpha(1_{F_\ell}\mu) \le U_\alpha(1_{F_\ell}\nu)$ holds for any $\alpha>0$ and $\ell$. 

\item[\((2)\)] {\sf (Convex cone)} If $\mu, \nu \in{\cal S}(\{F_\ell\})$ and $a, b \ge 0$, then $a \mu + b \nu \in{\cal S}(\{F_\ell\})$ and $U_\alpha(a1_{F_\ell}\mu+b1_{F_\ell}\nu) = a U_\alpha(1_{F_\ell}\mu) +bU_\alpha(1_{F_\ell}\nu)$ holds for $\alpha>0$ and 
$\ell$.

\item[\((3)\)] Let ${\cal B}_{b,+}(E)$ be the set of all nonnegative bounded Borel functions on $E$. Then ${\cal S}(\{F_{\ell}\})$ is closed under multiplication by
    ${\cal B}_{b,+}(E)$, namely, $f \mu$ belongs to ${\cal S}(\{F_\ell\})$ whenever $f\in {\cal B}_{b,+}(E)$ and  $\mu \in {\cal S}(\{F_\ell\})$.
\end{itemize}
Moreover the following  also holds:
\begin{itemize}
\item[\((4)\)] Let $\{G_\ell\}$ be another nest. Assume $\mu \in{\cal S}(\{F_\ell\})$ and $\nu \in {\cal S}(\{G_\ell\})$. Then $\mu+\nu \in {\cal S}(\{F_\ell\cap G_\ell\})$. 
\end{itemize}
\end{prop}

\begin{rem} We can replace ${\cal B}_{b,+}(E)$  in {\sf (3)} with 
the class of nonnegative Borel functions $f$ defined on $E$ that are bounded on each $F_\ell$.
\end{rem}

\begin{thm}[{\it c.f.} {\cite[Theorem 5.1]{NTTU25}}]
Let $\{F_\ell\}$ be a nest. Take $\{\mu_n, \mu\}  \subset \mathcal{S}(\{F_\ell\})$ 
and denote by $\{{\sf A}^n, {\sf A}\}$ their PCAFs. Assume that $\mu_n$ converges to $\mu$ 
in $\rho$ with respect to the nest $\{F_\ell\}$. Then there exists a subsequence 
$\{n_k\}$ such that for any $\ell \in \mathbb{N}$,
\[\mathbb{P}_x \Big( (1_{F_\ell} {\sf A}^{n_k})_t \to ( 1_{F_\ell} {\sf A} )_t \;\; \text{locally uniformly on } [0,\infty) \text{ as } n_k \to \infty \Big) = 1
\]
for q.e.\ $x \in E$.

If, in addition,
\[
\lim_{\ell \to \infty} \sup_{\nu \in S^1_{c,00},\, n \in \mathbb{N}}
\int_{F_\ell^c} U_1 \nu(x)\,\mu_n(dx) = 0
\quad \text{and} \quad
\lim_{\ell \to \infty} \sup_{\nu \in S^1_{c,00}}
\int_{F_\ell^c} U_1 \nu(x)\,\mu(dx) = 0,
\]
then there exists a further subsequence $\{n'_k\}$ of $\{n_k\}$ such that
\[\mathbb{P}_x \Big( {\sf A}^{n'_k}_t \to {\sf A}_t \;\; \text{locally uniformly on } [0,\infty) \text{ as } n'_k \to \infty \Big) = 1
\]
for q.e. $x \in E$.
Here ${\cal S}_{c, 00}^1:= \big\{ \mu \in {\cal S} : \,  {\sf supp}[\mu] \mbox{ is compact}, \ \mu(E)\le 1 \ 
\mbox{and} \ \|U_1\mu\|_\infty \le 1 \big\}$.
\end{thm}

\section{Convergence of Smooth Measures through the PCAFs Attached to a Nest} 

In this section, we utilize the continuity of the Revuz correspondence on the spaces attached to a nest to derive the convergence of smooth measures $\{\mu_n\}$ from the convergence of their associated PCAFs.

\begin{thm} \label{thm:main}
Let $\{\mu_n, \, \mu\}$ be smooth measures belonging to ${\cal S}(\{F_\ell\})$ and denote by 
$\{{\sf A}^n, \, {\sf A}\}$ the associated PCAFs.  
Assume the following conditions:

\begin{itemize}
\item[] \hspace*{-15pt} {\sf (A1)} for any $\ell \in{\mathbb N}$,
the functions  $\dis x\mapsto {\mathbb E}_x \big[(1_{F_\ell}{\sf A}^n)_t]$ converge to the function
 $x\mapsto  {\mathbb E}_x\big[ (1_{F_\ell}{\sf A})_t \big]$ in $L^1_{\sf loc}(E;m)$ 
 locally uniformly in $t$ on $[0, \infty)$ as $n\to \infty$ in the sense that, for any $t>0$ and 
 a compact set  $K\subset E$, 
$$
\lim_{n\to\infty} \sup_{0\le s\le t} \int_K \Big| {\mathbb E}_x\big[(1_{F_\ell}{\sf A}^n)_s]-
{\mathbb E}_x\big[ (1_{F_\ell}{\sf A})_s \big]\Big| m(dx)=0;
$$

\item[] \hspace*{-15pt} {\sf (A2)} for any $\alpha>0$ and $\ell \in {\mathbb N}$, the functions 
$\dis x\mapsto {\mathbb E}_x\Big[ \int_t^\infty \!\! e^{-\alpha s} d(1_{F_\ell}{\sf A}^n)_s\Big] 
\,\vee\, {\mathbb E}_x\Big[ \int_t^\infty \!\! e^{-\alpha s} d(1_{F_\ell}{\sf A})_s\Big]$ converge to 
$0$ in $L^1_{\sf loc}(E; m)$ {\it uniformly} in $n$  as $t\to\infty$ in the sense that, for any compact 
set $K\subset E$, 
$$
\lim_{t\to\infty} \sup_{n\in{\mathbb N}} \int_K \Big({\mathbb E}_x 
\Big[ \int_t^\infty e^{-\alpha s} d(1_{F_\ell}{\sf A}^n)_s\Big]  \vee 
{\mathbb E}_x\Big[ \int_t^\infty e^{-\alpha s} d(1_{F_\ell}{\sf A})_s\Big]\Big) m(dx) =0.
$$
\end{itemize}

\medskip
Then, $\{\mu_n\}$  converges to $\mu$ vaguely on the nest $\{F_\ell\}$ in the resolvent 
sense$:$ 
\begin{equation} \label{vcnr}
\lim_{n\to \infty} \int_{F_\ell} R_\alpha f(x)\mu_n(dx) = \int_{F_\ell} R_\alpha f(x) \mu(dx), \quad 
{\rm for \ every} \ \alpha>0, \ \ell \in{\mathbb N} \ {\rm and} \ f\in {\cal B}_{b, 0}(E), 
\end{equation}
where ${\cal B}_{b,0}(E):={\cal B}_{b,0}(E;m)$ is the set of all bounded Borel functions $f$ on $E$ 
such that ${\sf supp}[fm]$ is compact. 

Suppose, in addition, 
 \begin{equation} \label{1-pot}
M_\ell:=\sup_{x\in F_\ell} \Big(\sup_n \mathbb{E}_x 
\Big[ \int_0^\infty e^{-t} d\big(1_{F_\ell}{\sf A}^n\big)_t\Big]
\vee  \mathbb{E}_x\Big[ \int_0^\infty e^{-t} d\big(1_{F_\ell}{\sf A}\big)_t\Big] \Big) <\infty
 \end{equation}
holds for each $\ell\in{\mathbb N}$. Then {\it  $\{\mu_n\}$ converges weakly to $\mu$ 
on the nest $\{F_\ell\}$}. 
\end{thm}

\begin{rem}
\begin{enumerate}

\item[(1)]  We introduce  a notion ``vague convergence on the nest in the resolvent sense''  defined in (4.1) because smooth measures in our broad class are not necessarily Radon or finite measures, making the classical vague convergence inapplicable directly.

\item[(2)] The assumption that all measures $\{\mu_n\}$ and $\mu$ share the same nest $\{F_\ell\}$ might seem restrictive at first glance, but it is highly natural from the viewpoint of practical applications. In many applications (such as approximating highly singular objects like nowhere Radon measures), the target limit measure $\mu$ is given first, and a sequence of better-behaved approximating measures $\mu_n$ is constructed. In such scenarios, one can construct and fix a nest based on the limit measure $\mu$, and then evaluate the approximating sequence on this common nest (as we will demonstrate later in Section \ref{sec_example}). Thus, formulating the convergence on a fixed nest captures the very essence of such approximation procedures.

Furthermore, the assumptions (A1) and (A2) might appear somewhat technical, but they are in fact natural conditions arising from the intrinsic structure of the Revuz correspondence. Since our Dirichlet form and the Hunt process are defined on $L^2(E; m)$, the Revuz correspondence tightly links a smooth measure and its PCAF through the mediation of the underlying measure $m$ 
(as seen in Theorem \ref{thm-Revuz}). Therefore, to capture the vague convergence on the nest in the resolvent sense (which involves the resolvent $R_\alpha f$), it is most natural to evaluate the convergence of their corresponding PCAFs by ``lifting'' them to the integral with respect to $m$. Since the PCAF side involves a time integral over $[0, \infty)$, establishing its convergence naturally requires controlling both the local convergence on finite time intervals $[0, t]$ (which corresponds to (A1)) and the infinite tail $[t, \infty)$ (which corresponds to (A2)). Under this $m$-mediated framework, these conditions optimally bridge the gap between PCAFs and measures.

\end{enumerate}

\end{rem}

\begin{proof}[Proof of Theorem \ref{thm:main}]
Take $\ell \in{\mathbb N}$ and fix it. Let $f\in {\cal B}_{b,0}(E)$ and denote by
  $K:={\rm supp}[fm]$. Then, by virtue of \cite[Theorem 5.1.3 (iv)]{FOT11}, we see that, for any $\alpha>0,  t>0$ and $n\in{\mathbb N}$,   
\begin{align*}
& \!\!\!\!\! \Big| \int_{F_\ell}  R_\alpha f(x) \mu_n(dx) 
- \int_{F_\ell}  R_\alpha f(x) \mu(dx) \Big|  = \Big| \int_E R_\alpha f(x) (1_{F_\ell}\mu_n)(dx) 
- \int_E R_\alpha f(x) (1_{F_\ell}\mu)(dx) \Big| \\
& = \Big| \int_E f(x) {\mathbb E}_x\Big[ \int_0^\infty e^{-\alpha s} d(1_{F_\ell}{\sf A}^n)_s\Big] m(dx) -
\int_E f(x) {\mathbb E}_x\Big[ \int_0^\infty e^{-\alpha s} d(1_{F_\ell}{\sf A})_s\Big] m(dx) \Big| \\
& = \Big| \int_K f(x) {\mathbb E}_x\Big[ \int_0^\infty e^{-\alpha s} d(1_{F_\ell}{\sf A}^n)_s\Big] m(dx) -
\int_K  f(x) {\mathbb E}_x\Big[ \int_0^\infty e^{-\alpha s} d(1_{F_\ell}{\sf A})_s\Big] m(dx) \Big| \\
& \le  \bigg| \int_K f(x) \Big({\mathbb E}_x\Big[ \int_0^t e^{-\alpha s} d(1_{F_\ell}{\sf A}^n)_s\Big]
-{\mathbb E}_x\Big[ \int_0^t e^{-\alpha s} d(1_{F_\ell}{\sf A})_s\Big] \Big) m(dx) \biggl|  \\
& \qquad + \biggl|\int_K f(x) {\mathbb E}_x\Big[ \int_t^\infty e^{-\alpha s} 
d(1_{F_\ell}{\sf A}^n)_s\Big] m(dx)  \biggl| + \biggl| \int_K f(x) {\mathbb E}_x 
\Big[ \int_t^\infty e^{-\alpha s} d(1_{F_\ell}{\sf A})_s\Big] m(dx)  \biggl|. 
\end{align*}
By using the integration by parts with respect to the Stieltjes integrals, the right hand 
side of the above is dominated by 
\begin{align*}
& \biggl| \int_K  f(x) \Big( 
{\mathbb E}_x\Big[ e^{-\alpha t} {(1_{F_\ell}{\sf A}^n})_t  
+\alpha \int_0^t e^{-\alpha s} (1_{F_\ell}{\sf A}^n)_sds\Big] \\
& \qquad  -{\mathbb E}_x\Big[ e^{-\alpha t} {(1_{F_\ell}{\sf A}})_t  
+\alpha \int_0^t e^{-\alpha s} (1_{F_\ell}{\sf A})_sds\Big] \Big)m(dx) \biggl| \\
& \qquad  + 2 \|f\|_\infty   \sup_{n\in{\mathbb N}}  
\int_K \Big( {\mathbb E}_x\Big[ \int_t^\infty e^{-\alpha s} d(1_{F_\ell}{\sf A}^n)_s\Big] \vee 
{\mathbb E}_x\Big[ \int_t^\infty e^{-\alpha s} d(1_{F_\ell}{\sf A})_s\Big]\Big) m(dx) \\
& \le  e^{-\alpha t} \|f\|_\infty \int_K \big|  {\mathbb E}_x\big[ (1_{F_\ell}{\sf A}^n)_t \big] 
-{\mathbb E}_x\big[(1_{F_\ell}{\sf A})_t\big] \big| m(dx)   \\ 
& \qquad  + \alpha  \|f\|_\infty \int_0^t e^{-\alpha s} \Big( \int_K \big| {\mathbb E}_x\big[ (1_{F_\ell}{\sf A}^n)_s \big] - 
 {\mathbb E}_x\big[ (1_{F_\ell}{\sf A})_s \big]\big| m(dx) \Big) ds \\
 & \qquad + 2 \|f\|_\infty  \sup_{n\in{\mathbb N}} \int_K 
\left(  {\mathbb E}_x\Big[ \int_t^\infty e^{-\alpha s} d(1_{F_\ell}{\sf A}^n)_s\Big]
 \vee {\mathbb E}_x\Big[ \int_t^\infty e^{-\alpha s} d(1_{F_\ell}{\sf A})_s\Big] \right) m(dx) \\
& = 2 \|f\|_\infty \sup_{0\le s\le t}  \int_K \big| {\mathbb E}_x\big[(1_{F_\ell}{\sf A}^n)_s\big] - 
{\mathbb E}_x\big[(1_{F_\ell}{\sf A})_s \big] \big| m(dx)   \\
& \qquad + 2 \|f\|_\infty  \sup_{n\in{\mathbb N}} \int_K 
\left(  {\mathbb E}_x\Big[ \int_t^\infty e^{-\alpha s} d(1_{F_\ell}{\sf A}^n)_s\Big]
 \vee {\mathbb E}_x\Big[ \int_t^\infty e^{-\alpha s} d(1_{F_\ell}{\sf A})_s\Big] \right) m(dx).
\end{align*}
Then by virtue of the conditions {\sf (A1)} and {\sf (A2)},  tending $n\to \infty$ 
 and then $t\to \infty$,  it follows that the right hand side goes to $0$. Thus 
$$
\lim_{n\to \infty} \int_{F_\ell} R_\alpha f(x)\mu_n(dx) = \int_{F_\ell} R_\alpha f(x) \mu(dx), \ \  {\rm for  \ every \  }  
\alpha>0, \ \ell \in{\mathbb N} \ {\rm and} \ f \in {\cal B}_{b,0}(E).
$$

We now turn to the proof of the weak convergence of \(\mu_n\) to \(\mu\) on the nest \(\{F_\ell\}\) under the additional assumption that \eqref{1-pot} holds. 
We assume  that \eqref{1-pot} holds for each $\ell \in{\mathbb N}$. We will show that the $1$-potentials $\{U_1(1_{F_\ell}\mu_n)\}$ are bounded with respect  to $\form_1$-norm in $\dom$.
For this, we note that since $F_\ell$ is compact and the form $(\form, \dom)$ is regular,  
 there exists a function $\varphi\in \dom \cap C_0(E)$ such that $0\le \varphi \le 1$ and $\varphi =1$ on $F_\ell$.
The fact that $\mathbb{E}_{\cdot}\Big[ \int_0^\infty e^{-t}d(1_{F_\ell}{\sf A}^n)_t\Big]$ 
is a quasi continuous modification of the $1$-potential  $U_1(1_{F_\ell}\mu_n)$ of the measure $1_{F_\ell}\mu_n$ 
for each $\ell$ and $n$, combined with the characterization of the $1$-potentials, yields:
\begin{align*}
\form_1(U_1(1_{F_\ell}\mu_n), U_1(1_{F_\ell}\mu_n))
 & =\int_{F_\ell} {
 {U_1(1_{F_\ell}\mu_n)}(x)} \mu_n(dx) = \int_{F_\ell} \mathbb{E}_x\Big[\int_0^\infty e^{-t}d(1_{F_\ell}{\sf A}^n)_t\Big] \mu_n(dx)  \\
& \le M_\ell \int_{F_\ell} \varphi(x) \mu_n(dx) =M_\ell\, \form_1(\varphi,  U_1(1_{F_\ell}\mu_n)) \\
& \le  M_\ell \sqrt{\form_1(\varphi, \varphi)} \sqrt{\form_1(U_1(1_{F_\ell}\mu_n), U_1(1_{F_\ell}\mu_n))}.
\end{align*}
By dividing both sides of the last inequality by $\sqrt{\form_1(U_1(1_{F_\ell}\mu_n), U_1(1_{F_\ell}\mu_n))}$, we get
\begin{align*}
 \|U_1(1_{F_{\ell}} \mu_n)\|_{\form_1} \le M_{\ell} \|\varphi\|_{\form_1}.
\end{align*}
This implies the boundedness of the sequence:
$$
\sup_n \form_1(U_1(1_{F_\ell}\mu_n), U_1(1_{F_\ell}\mu_n)) \le M_\ell^2 \form_1(\varphi, \varphi) <\infty.
$$
A similar argument for the limit potential $U_1(1_{F_\ell}\mu)$ establishes its boundedness. Hence, we conclude that both the sequence of potentials and its limit are bounded in the $\form_1$-norm.

Next, we show that the $1$-potentials $\{U_1(1_{F_\ell}\mu_n)\}$ converge weakly to the $1$-potential of $1_{F_\ell}\mu$ 
in $\dom$ with respect to $\form_1$-norm.  This weak convergence, along with the boundedness we have already shown, 
 allows us to conclude the weak convergence of the measures $\{\mu_n\}$ on the nest $\{F_\ell\}$ by repeating the proof of
 Proposition \ref{prop:weak}.   
 
To prove this, take any $\varphi \in \dom$ and $\ell \in\mathbb{N}$. Since the space $R_1(\dom\cap C_0(E))$ is dense 
in $\dom$  with respect to $\form_1$-norm, for any given $\vareps>0$  and $\ell \in{\mathbb N}$, 
 there exists a $g\in \dom\cap C_0(E) \, \big( \! \subset {\cal B}_{b,0}(E)\big)$ 
 such that $\form_1(\varphi-R_1g, \varphi-R_1g)<\vareps^2/4 M_\ell$.  
We evaluate the difference $|\form_1(\varphi, U_1(1_{F_\ell}\mu_n)) - \form_1(\varphi, U_1(1_{F_\ell}\mu))|$. 
By the triangle inequality, we can bound this difference by three terms:
\begin{align*}
\lefteqn{\left| \form_1(\varphi, U_1(1_{F_\ell}\mu_n)) - \form_1(\varphi, U_1(1_{F_\ell}\mu))\right|  \le \left| \form_1(\varphi -R_1g, U_1(1_{F_\ell}\mu_n))\right|}\\
& + \left| \form_1(R_1g, U_1(1_{F_\ell}\mu_n)) - \form_1(R_1g, U_1(1_{F_\ell}\mu)) \right| + \left| \form_1(R_1g - \varphi, U_1(1_{F_\ell}\mu))\right).
\end{align*}
Applying the Schwarz inequality to the first and third terms, and utilizing the boundedness of the potentials established in the previous step, we can estimate this expression.  Then the right hand side is bounded by
\begin{align*}
& \sqrt{\form_1(\varphi -R_1g)} \sqrt{ \form_1(U_1(1_{F_\ell}\mu_n))} + \left| \form_1(R_1g, U_1(1_{F_\ell}\mu_n)) - \form_1(R_1g, U_1(1_{F_\ell}\mu)) \right| \\
& \qquad + \sqrt{\form_1(R_1g - \varphi)} \sqrt{ \form_1(U_1(1_{F_\ell}\mu))} \\
 & \le  2\varepsilon M_\ell
 + \left| \int_{F_\ell} R_1 g d\mu_n - \int_{F_\ell} R_1 g d\mu \right|.
\end{align*}
Since the measures $\{\mu_n\}$ converge to $\mu$ vaguely
on the nest $\{F_\ell\}$ in the resolvent sense (see \eqref{vcnr}), 
the last term on the right-hand side goes to zero as $n\to\infty$. As $\vareps$ is an arbitrary positive number, 
this implies that the left-hand side must converge to zero.  
Therefore, the sequence of potentials $\{U_1(1_{F_\ell}\mu_n)\}$ converges to $U_1(1_{F_\ell}\mu)$ 
{\it weakly in} $\dom$ with respect to $\form_1$-norm as $n\to\infty$. 
\end{proof}

The following corollary immediately follows from Theorem \ref{thm:main}  and Proposition \ref{prop:weak}.

\begin{cor} \label{cor:1}
Let $\{\mu_n, \mu\} \subset {\cal S}(\{F_\ell\})$ for  an exhaustion $\{F_\ell\}$ of $E$. Assume that the conditions {\sf (A1), (A2)} and \eqref{1-pot} hold for $\{F_\ell\}$.  
Then the measures  $\{\mu_n\}$ converge to $\mu$ {\it vaguely}.
\end{cor}

The next corollary shows that the vague convergence holds for measures of finite energy integrals, given 
the boundedness of their 1-potentials in the $\form_1$-norm and the conditions of Theorem \ref{thm:main}:  

\begin{cor} \label{cor:2}
 Let $\{F_\ell\}$ be an exhaustion of $E$, $\{\mu_n\} \subset {\cal S}_0$ and $\mu\in \mathcal{S}(\{F_\ell\})$. Assume that the conditions 
{\sf (A1)}, {\sf (A2)} and the following condition 
 \begin{equation} \label{1-bdd}
M:= \sup_n \form_1(U_1\mu_n,U_1\mu_n)<\infty
\end{equation}
holds, then the measures $\{\mu_n\}$ converge to $\mu$ {\it vaguely}, and the measure $\mu$ also belongs to ${\cal S}_0$.
Moreover, the potentials $\{U_1\mu_n\}$ converge to $U_1\mu$ {\it weakly in} $(\dom, \form_1)$: 
$$
 \lim_{n\to \infty} \form_1(\varphi, U_1\mu_n) = \form_1(\varphi, U_1\mu),
 \qquad \text{for any }\varphi\in \dom.
$$
\end{cor}
\begin{proof}
We show $\mu \in {\cal S}_0$, since the vague convergence of $\{\mu_n\}$ 
and the weak convergence of $\{U_1\mu_n\}$ are shown in a way  similar to the previous one by replacing 
$M_\ell$ with $M \vee \form_1(U_1(1_{F_\ell}\mu),U_1(1_{F_\ell}\mu))$.  

Noting that $1_{F_\ell}  \mu_n$ and $\mu_n$ are all belonging  to ${\cal S}_0$ for any $\ell$ and $n$,  we have 
\begin{align*}
\Big|  \int_{F_\ell} R_1f(x) \mu_n(dx) \Big|  & = \big| \form_1(R_1 f, U_1(1_{F_\ell} \mu_n))\big| 
 \le \sqrt{\form_1(R_1f, R_1f)} \sqrt{ \form_1(U_1(1_{F_\ell} \mu_n), U_1(1_{F_\ell} \mu_n))}  \\
 & \le \sqrt{\form_1(R_1f, R_1f)} \sqrt{ \form_1(U_1\mu_n, U_1\mu_n)}
 \leq {
 \sqrt{M}} \sqrt{\form_1(R_1f, R_1f)}
\end{align*}
for each $f\in \dom\cap C_0(E) \, \big(\subset {\cal B}_{b,0}(E)\big)$. 
Here we used the following inequalities holding  for $\nu \in{\cal S}_0$:
$$
U_1(1_{F_\ell} \nu) \le U_1\nu
$$
and
$$
\form_1(U_1(1_{F_\ell} \nu), U_1(1_{F_\ell} \nu)) \le \form_1(U_1\nu, U_1\nu)
$$
(see  Proposition \ref{smoooth_alg} (1)).  Since $\{\mu_n\}$ converges to $\mu$ vaguely on 
the nest $\{F_\ell\}$ in the resolvent sense (Theorem \ref{thm:main}), we have 
$$
 \Big|\int_E R_1f(x) \mu(dx)\Big|  \le {
 \sqrt{M}} \sqrt{\form_1(R_1f, R_1f)},  \quad f\in \dom\cap C_0(E).
$$
Since $R_1(\dom\cap C_0(E))$ is dense in $\dom$, it follows that $\mu$ is a measure of finite energy integrals.

For any \(\varphi \in \mathcal{F}\) and \(\varepsilon>0\), by the regularity, we take \(\varphi_\varepsilon \in \mathcal{F}\cap C_0(E)\) satisfying  \(\sqrt{\mathcal{E}_1(\varphi-\varphi_\varepsilon,\varphi-\varphi_\varepsilon)}\leq \varepsilon\). Then  we have
\begin{eqnarray*}
| \form_1(\varphi, U_1\mu_n) -\form_1(\varphi, U_1\mu)|& \leq &| \form_1(\varphi-\varphi_\varepsilon, U_1\mu_n) | + | \form_1(\varphi-\varphi_\varepsilon, U_1\mu) | + | \form_1(\varphi_\varepsilon, U_1\mu_n-U_1\mu)|\\
&\leq & 2\sqrt{M} \sqrt{\mathcal{E}_1(\varphi-\varphi_\varepsilon,\varphi-\varphi_\varepsilon)} + \left|\int_E \varphi_\varepsilon d\mu_n - \int_E \varphi_\varepsilon d\mu\right|\\
&\leq & 2\sqrt{M} \varepsilon + \left|\int_E \varphi_\varepsilon d\mu_n - \int_E \varphi_\varepsilon d\mu\right|,
\end{eqnarray*}
and so, by letting \(n\to \infty\) and then \(\varepsilon \searrow 0\), the potentials $\{U_1\mu_n\}$ converge to $U_1\mu$ weakly in $(\dom, \form_1)$. 
\end{proof}

We consider another perspective of conditions in Theorem \ref{thm:main}. 
 \begin{prop} \label{propforA2}
For \(\mu_n, \mu \in \mathcal{S}(\{F_\ell\})\) satisfying \eqref{1-pot} for each \(\ell \in \mathbb{N}\), \(\sf (A2)\) holds. 
\end{prop}
 \begin{proof}
 For \(\mu_n, \mu \in \mathcal{S}(\{F_\ell\})\) satisfying \eqref{1-pot},
 we have already proved in the proof of Theorem \ref{thm:main} that
 \begin{align}
   \sup_{n \in \bN} \cE_1(U_1(1_{F_{\ell}}\mu_n), U_1(1_{F_{\ell}}\mu_n))
  \vee \cE_1(U_1(1_{F_{\ell}} \mu), U_1(1_{F_{\ell}} \mu)) < \infty.
 \end{align}
 By \cite[Remark 3.2]{NTTU25}, it holds that
 \begin{align}
    M_{\ell} :=\sup_{n}\mathcal{E}_{\alpha}
  (U_{\alpha}(1_{F_\ell} \mu_n),U_{\alpha}(1_{F_\ell} \mu_n))
  \vee \mathcal{E}_{\alpha}(U_{\alpha}(1_{F_\ell}\mu),
  U_{\alpha}(1_{F_\ell}\mu)) <\infty\quad
  \text{for any}\ \alpha > 0.
 \end{align}
Let \({\sf A}^n, {\sf A}\) be PCAFs corresponding to \(\mu_n, \mu\), respectively. By the Markov property and symmetry of $p_t$, for a compact set \(K\), we have 
 \begin{align*}
  \int_K \mathbb{E}_x\left[\int_t^{\infty} e^{-\alpha s} d(1_{F_{\ell}}{\sf A}^n)_s \right]
  dm(x) 
    &= \int_{K} \bE_x\left[\int_{0}^{\infty} e^{-\alpha (s+t)} d(1_{F_{\ell}}{\sf A}^n)_{s+t}\right]dm(x)\\
  &= e^{-\alpha t} \int_{K} \bE_x\left[\int_{0}^{\infty}
				       e^{-\alpha s} d (1_{F_{\ell}}{\sf A}^n_t+1_{F_{\ell}} {\sf A}_s^n \circ
  \theta_t)\right] m(dx) \\
  &= e^{-\alpha t} \int_E 1_K \cdot p_t
  \left(\mathbb{E}_{\cdot}\left[
  \int_{0}^{\infty} e^{-\alpha s} 1_{F_{\ell}}(X_s) d{\sf A}_s^n
  \right]\right)(x)m(dx) \\
  &= e^{-\alpha t} \int_E p_t1_K \cdot 
  \mathbb{E}_{x}\left[
  \int_{0}^{\infty} e^{-\alpha s} 1_{F_{\ell}}(X_s) d{\sf A}_s^n
  \right]m(dx) \\
  &= e^{-\alpha t} \mathcal{E}_{\alpha}(U_{\alpha}(1_{F_{\ell}} \mu_n),
  R_{\alpha}p_t1_K)\\
  &\le e^{-\alpha t} \sqrt{\mathcal{E}_{\alpha}(U_{\alpha}(1_{F_l
  }\mu_n))}
  \sqrt{\mathcal{E}_{\alpha}(R_{\alpha}p_t1_K)}.
 \end{align*}
 Now, using the symmetry of $p_t$ again and the Markovian property of $p_t$
 and $R_{\alpha}$, we have
 \begin{align}
  \mathcal{E}_{\alpha}(R_{\alpha}p_t1_K,R_{\alpha}p_t1_K) &=
  \int_E p_t1_K \cdot R_{\alpha}p_t1_K m(dx)
  = \int_K p_t R_{\alpha}p_t1_K m(dx)\\
  &\le  \int_K \frac{1}{\alpha} m(dx) = \frac{m(K)}{\alpha} < \infty.
 \end{align}
So we have
 \[
 \sup_n \int_K \mathbb{E}_x\left[\int_t^{\infty} e^{-\alpha s} d{\sf A}_s^n
 \right] m(dx) 
 \leq e^{-\alpha t} \cdot \frac{\sqrt{M_{\ell}} \cdot  m(K)}{\alpha}
 \]
  and this converges to \(0\) as \(t\to \infty.\) Similarly, \(\sup_n \int_K \mathbb{E}_x\left[\int_t^{\infty} e^{-\alpha s} d{\sf A}_s^n
 \right] m(dx) \) converges to  \(0\) as \(t\to \infty.\)
 \end{proof}

In our previous work \cite{NTTU25}, we introduced a class of functions 
$\dot{L}^1_{\sf loc}(E;m)$, which are locally in $L^1(E;m)$ in the broad sense. 
We now extend this concept to a similar class for functions that are locally in 
 $L^p(E;m)$ in the broad sense for  $1<p \le \infty$. 

A measurable function $f$ defined {\it q.e.} on $E$  is said to be {\it locally 
in $L^p(E;m)$ in the broad sense} (denoted as $f\in \dot{L}^p_{\sf loc}(E;m)$) 
if there exists a nest $\{F_\ell\}_\ell$ such that 
$f|_{F_\ell}\in L^p(F_\ell;m)$ for each $\ell$.  

For a nest $\{F_\ell\}$, we also use the notation $f\in L^p(\{F_\ell\};m)$ to indicate that 
$f|_{F_\ell}\in L^p(F_\ell;m)$ for each $\ell$.  By definition, it immediately follows that $\dot{L}^p_{\sf loc}(E;m) = \bigcup_{\{F_\ell\} \in \nest} L^p(\{F_\ell\};m)$.  We then observe the following inclusions: for 
$1\le p \le q \le \infty$, 
$$
L^q_{\sf loc}(E;m)\subset \dot{L}^q_{\sf loc}(E;m) 
= \bigcup_{\{F_\ell\} \in \nest} 
L^q(\{F_\ell\};m) \subset \dot{L}^p_{\sf loc}(E;m)  \subset \dot{L}^1_{\sf loc}(E;m).
$$
For \(p\geq 2\), we may identify $L^p(\{F_\ell\};m )\cap{\cal B}_+(E)$ with a class of measures in \(\mathcal{S}(\{F_{\ell}\})\). Indeed, for any \(f \in L^p(\{F_\ell\};m )\cap{\cal B}_+(E)\) and \(g\in \mathcal{F}\), we have 
$$
\int_{F_{\ell}} |g|f dm \leq \|f\|_{L^2(F_{\ell})} \|g\|_{L^2(F_{\ell})} \leq \|f\|_{L^p(F_{\ell})} m(F_{\ell})^{1/p^*} \sqrt{\mathcal{E}_1(g,g)}
$$
where \(p^*=2p/(p-1)\), and so \(f dm\in \mathcal{S}(\{F_{\ell}\}) \).

For \(1\leq p< 2\), \(L^p(\{F_\ell\};m )\cap{\cal B}_+(E)\) cannot be identified with the class of measures in \(\mathcal{S}(\{F_{\ell}\})\). For example, for \(d\geq 3\) and \((\mathcal{E}, \mathcal{F})=(\frac{1}{2}\mathbb{D}, H^1(\mathbb{R}^d))\), let \(f(x):=|x|^{-\beta}\) for \(\frac{d+2}{2}<\beta <d\) and \(F_{\ell}:=\{|x|\leq \ell\}\). 
Then, for \(p<d/\beta\), \(f \in L^p(\{F_\ell\};dx)\cap{\cal B}_+(\real^d)\) but \(fdx \not \in \mathcal{S}(\{F_{\ell}\})\) since, for 
$\varphi \in C_0(\real^d)$ with $\varphi \ge 0$  and \(\varphi_{\varepsilon}(x):=\varphi(x/\varepsilon), \ \vareps>0\), there exists \(C>0\) independent of \(\varepsilon \) such that 
\begin{eqnarray*}
\varliminf_{\varepsilon \to 0} \frac{\int_{F_{\ell}} \varphi_{\varepsilon}(x)f(x)dx}{\sqrt{\mathcal{E}_1(\varphi_{\varepsilon},\varphi_{\varepsilon})}} \geq \varliminf_{\varepsilon \to 0} \frac{\varepsilon^{d-\beta}}{\varepsilon^{(d-2)/2}} C =\varliminf_{\varepsilon \to 0}\varepsilon^{\frac{d+2}{2}-\beta} C =\infty.
\end{eqnarray*}
However,  according to \cite[Lemma 6.1]{NTTU25}, for a nonnegative Borel function $f$ on $E$ 
that belongs to $\dot{L}^p_{\sf loc}(E;m)$ for some $1\le p \le \infty$,  
the measure \(\mu:=fm\) is smooth and its associated PCAF is given by 
$$
{\sf A}_t:=\int_0^t f(X_s)ds.
$$
This means that for \(1 \leq p <2\) and \(f \in L^p(\{F_{\ell}\};m)\cap{\cal B}_+(E)\), \(\mu=f dm \) is smooth, but it may be attached to a different nest than \(\{F_{\ell}\}\). In particular, in the above example with \(f(x)=|x|^{-\beta}\) for \(\frac{d+2}{2}<\beta <d\), we have \(fdx \in \mathcal{S}(\{\tilde{F}_{\ell}\})\) where \(\tilde{F}_{\ell}:=\{1/\ell \leq |x| \leq \ell\}\). See also Example \ref{ex-01}.

\begin{prop}\label{prop:broad}
Let \(\{F_\ell\}\) be a compact nest and take nonnegative Borel functions 
$f, f_n \in L^p(\{F_\ell\};m)$  for some \(1 \leq p\le \infty\)  such that $\{f_n\}$ converges to $f$ 
{\it strongly} in $L^p(F_\ell;m)$ for each $\ell\in{\mathbb N}$.  For \(1\leq p<2,\) assume that there exist
an increasing sequence of positive numbers $\{a_\ell\}_{\ell=1}^\infty$ and a compact nest \(\{C_{\ell}\}_\ell\) such that
\begin{equation} \label{eq:nest}
C_{\ell} \subset \left\{ x \in   F_\ell :  \, \big(\sup_n f_n(x)\big) \vee f(x) \le a_\ell \right\} \quad {\rm for\ any\ } \ell \in{\mathbb N},
\end{equation}
and,  for \(2\leq p\leq \infty\), we set \(C_{\ell}:=F_{\ell}.\)  Then the conditions 
{\sf (A1)}, {\sf (A2)} and \eqref{1-pot} in Theorem \ref{thm:main} hold for the nest $\{C_{\ell}\}$, 
smooth measures $\mu:=fm$ and $\mu_n:=f_nm$, and their associated PCAFs \(\sf A^n, \sf A\). 
Consequently, $\mu_n$ {\it converges weakly to}  $\mu$ {\it on the nest} $\{C_{\ell}\}$.
\end{prop}

\begin{proof}  We set \(p^*:=p/(p-1)\) for \(1<p<\infty \), \(p^*:=\infty\) and \(1/p^*:=0\) for \(p=1\), 
and \(p^*:=1\) for \(p=\infty\). We first consider the case of \(1\leq p< 2\).   We note that the inequality 
\begin{align*}
\int_{C_{\ell}} f(x)^2m(dx) \vee \int_{C_{\ell}} f_n(x)^2m(dx) & \le a_\ell \Big( 
\int_{C_{\ell}} |f(x)|m(dx) \vee \int_{C_{\ell}} |f_n(x)|m(dx) \Big)  \\
& \leq  a_\ell m(C_{\ell})^{1/p^*} \big( \|f\|_{L^p(C_{\ell})} \vee \|f_n\|_{L^p(C_{\ell})} \big) <\infty
\end{align*}
 holds for  each $\ell$ and $n$, 
 and so $1_{C_{\ell}}\mu,\, 1_{C_{\ell}}\mu_n \in {\cal S}_0$. 
Thus this implies that $\{\mu_n, \mu\} \subset {\cal S}(\{C_{\ell}\})$

For each $\ell$ and $n$, by the definition of $C_{\ell}$, we have
$$
{\mathbb E}_x\Big[\int_0^\infty e^{-t}d(1_{C_{\ell}}{\sf A}^n)_t\Big] =
{\mathbb E}_x\Big[\int_0^\infty e^{-t}1_{C_{\ell}}(X_t)f_n(X_t)dt\Big]
\le a_\ell {\mathbb E}_x\Big[\int_0^\infty e^{-t} dt \Big] =a_\ell. 
$$
and this implies that \eqref{1-pot} holds for $\{C_{\ell}\}$.  Moreover  for any compact set $K\subset E$, $\ell \in{\mathbb N}$ 
and $t>0$, 
\begin{align*}
& \!\!\! \sup_{0\le s\le t} \int_K  \big|  {\mathbb E}_x\big[ (1_{C_{\ell}}{\sf A}^n)_s\big]  -
{\mathbb E}_x\big[ (1_{C_{\ell}}{\sf A})_s \big]  \big| m(dx) \\
& \le  \int_K  \sup_{0\le s\le t} {\mathbb E}_x\Big[ \Big| \int_0^s 1_{C_{\ell}}(X_u)f_n(X_u)du -
\int_0^s 1_{C_{\ell}}(X_u)f(X_u)du \Big| \Big] m(dx) \\
& \le \int_K  {\mathbb E}_x\Big[ \int_0^t  1_{C_{\ell}}(X_u)\big| f_n(X_u)-f(X_u)\big| du \Big] m(dx) \\
&  = \int_0^t \int_K p_u(1_{C_{\ell}}|f_n-f|)(x) m(dx) du  = \int_0^t \int_{C_{\ell}} |f_n(x)-f(x)|p_u1_K(x) m(dx) du \\ 
& \le t  \|f_n-f\|_{L^1(C_{\ell})} \leq  t  \|f_n-f\|_{L^p(C_{\ell})} m(C_{\ell})^{1/p^*}
\to 0 \quad  {\rm as} \ n\to \infty.
\end{align*}
Thus we see {\sf (A1)} and, by Proposition \ref{propforA2}, {\sf (A2)} hold for the nest $\{C_{\ell}\}$. 
Hence Theorem \ref{thm:main} implies that $\{\mu_n\}$ converges weakly to $\mu$ on the nest $\{C_{\ell}\}$.

For the case of \(2\leq p\leq \infty\), since \(L^p(\{F_{\ell}\};m)\) can be identified with the class of 
\(\mathcal{S}(\{F_{\ell}\})\), we have $\{\mu_n, \mu\} \subset {\cal S}(\{F_{\ell}\})$. 
In the same way as above, we can check that \eqref{1-pot}, {\sf (A1)} and {\sf (A2)} hold for the nest $\{F_{\ell}\}$ 
and so $\{\mu_n\}$ converges weakly to $\mu$ on the nest  $\{F_{\ell}\}$ for the case of 
\(2\leq p\leq \infty\). 
\end{proof}

\begin{rem}
\((1)\) Under the assumption in Proposition \ref{prop:broad}, we also see that \(\mu_n\) converges to \(\mu\) in \(\rho\) for the nest $\{C_{\ell}\}$.  Indeed, for \(1\leq p<2\) we have
\begin{align*}
\rho_0(1_{C_\ell} \mu_n, 1_{C_\ell} \mu)^2=\int_{C_{\ell}}(f_n-f)R_1(1_{C_\ell} f_n-1_{C_\ell} f)dm \leq  2a_{\ell} \|f_n-f\|_{L^1(F_{\ell})}\to 0
\end{align*}
and for \(2\leq p\leq \infty\),
\begin{align*}
\rho_0(1_{F_\ell} \mu_n, 1_{F_\ell} \mu)^2=\int_{F_{\ell}}(f_n-f)R_1(1_{F_\ell} f_n-1_{F_\ell} f)dm \leq  \|f_n-f\|_{L^2(F_{\ell})}^2 \to 0.
\end{align*} 
\((2)\) In the case where non-negative Borel functions $f_n$ converge to $f$ in $L^p(E;m)$ for some $2\le p<\infty$, we can directly verify that the conditions (A1), (A2) and \eqref{1-pot} in Corollary \ref{cor:1} are satisfied for any increasing sequence of compact sets $\{F_\ell\}$ that forms an exhaustion of $E$ (i.e., $\bigcup_{\ell=1}^\infty F_\ell = E$).  Since our underlying space $E$ is a locally compact separable metric space, it is then $\sigma$-compact, which guarantees that such an exhaustion $\{F_\ell\}$ always exists. This implies that $\mu_n=f_nm$ converges to $\mu=fm$ vaguely. Of course, we can also demonstrate this vague convergence directly even without invoking the nest $\{F_{\ell}\}$.
\end{rem}

\medskip
We close this section by proving the following theorem:

\begin{thm}\label{conv_rho_nest}
Let $\{F_\ell\}$ be a nest. Let $\{\mu_n\}$ and $\mu$ be smooth measures in ${\cal S}(\{F_\ell\})$, and let $\{{\sf A}^n\}$ and ${\sf A}$ be the corresponding PCAFs.  Assume that for each $\ell \in{\mathbb N}$, the following conditions hold:
\begin{itemize}
\item[\((1)\)] {\rm (Local convergence of 1-potentials)} For any $t>0$, 
$$\dis 
\lim_{n\to \infty} \sup_{x\in F_\ell} \Big|
{\mathbb E}_x \Big[ \int_0^t e^{-s} d(1_{F_\ell}{\sf A}^n)_s\Big]-
{\mathbb E}_x \Big[ \int_0^t e^{-s} d(1_{F_\ell}{\sf A})_s\Big]\Big|=0.
$$

\item[\((2)\)] {\rm (Uniform tail estimate)} \[\dis 
\lim_{t\to \infty} \sup_{x\in F_\ell,\, n \in{\mathbb N}} 
\Big( {\mathbb E}_x \Big[ \int_t^\infty e^{-s} d(1_{F_\ell}{\sf A}^n)_s \Big]  \vee
 {\mathbb E}_x \Big[ \int_t^\infty e^{-s}d(1_{F_\ell}{\sf A})_s \Big] \Big) =0.
\]
\end{itemize}
Then $\mu_n$ converges to $\mu$ in the metric $\rho$ with respect to $\{F_\ell\}_\ell\).
\end{thm}

\begin{proof} 
We first observe that $\sup_n \mu_n(F_\ell) < \infty$ holds for each $\ell$. We fix \(\ell \in \mathbb{N}\). Since $F_\ell$ is compact and $(\mathcal{E}, \mathcal{F})$ is regular, we take $\phi \in \mathcal{F} \cap C_0(E)$ satisfying $\phi = 1$ on $F_\ell$ and $0 \le \phi \le 1$ on $E$. By assumption (2), for any given $\varepsilon > 0$, we can choose a sufficiently large $t > 0$ such that 
$\sup_{x\in F_\ell, n\in\mathbb{N}} {\mathbb E}_x[\int_t^\infty e^{-s}d(1_{F_\ell}{\sf A}^{n})_s] \le \varepsilon/2$. For this fixed $t$, 
assumption (1) implies that for sufficiently large $n$, 
$\sup_{x\in F_\ell} |{\mathbb E}_x[\int_0^t e^{-s}d(1_{F_\ell}{\sf A}^{n})_s] - {\mathbb E}_x[\int_0^t e^{-s}d(1_{F_\ell}{\sf A})_s]| \le \varepsilon/2$. 
Consequently, for such $n$ and q.e. $x \in F_\ell$,  recalling that ${\mathbb E}_\cdot[\int_0^\infty e^{-s}d(1_{F_\ell}{\sf A}^{n})_s]$ and ${\mathbb E}_\cdot[\int_0^\infty e^{-s}d(1_{F_\ell}{\sf A})_s]$ are quasi-continuous versions of the 1-potentials $U_1(1_{F_\ell}\mu_n)$ and $U_1(1_{F_\ell}\mu)$, respectively, we have 
$$
U_1(1_{F_\ell}\mu_n)(x)  \le {\mathbb E}_x\left[ \int_0^t e^{-s}d(1_{F_\ell}{\sf A})_s \right] + \varepsilon \le U_1(1_{F_\ell}\mu)(x)  + \varepsilon
$$
for q.e.\(x\). Integrating this inequality over $F_\ell$ with respect to $\mu_n$, we obtain
\begin{align}
\nonumber \lefteqn{\mathcal{E}_1(U_1(1_{F_\ell}\mu_n), U_1(1_{F_\ell}\mu_n)) = \int_{F_\ell} U_1(1_{F_\ell}\mu_n)(x)  \mu_n(dx)} \\
\nonumber & \le \int_{F_\ell} U_1(1_{F_\ell}\mu)(x)  \mu_n(dx) 
+ \varepsilon \mu_n(F_\ell) \ =\  \mathcal{E}_1(U_1(1_{F_\ell}\mu), U_1(1_{F_\ell}\mu_n)) + \varepsilon \mu_n(F_\ell) \\
\nonumber &\le \sqrt{\mathcal{E}_1(U_1(1_{F_\ell}\mu), U_1(1_{F_\ell}\mu))} \sqrt{\mathcal{E}_1(U_1(1_{F_\ell}\mu_n), U_1(1_{F_\ell}\mu_n))}  \\
& \qquad + \varepsilon \sqrt{\mathcal{E}_1(\phi, \phi)} \sqrt{\mathcal{E}_1(U_1(1_{F_\ell}\mu_n), U_1(1_{F_\ell}\mu_n))}, \label{eq:4.8-1}
\end{align}
where we used 
\begin{equation}
\mu_n(F_\ell) \le \int_{F_\ell} \phi d\mu_n = \mathcal{E}_1(U_1(1_{F_\ell}\mu_n), \phi) 
\le \sqrt{\mathcal{E}_1(\phi, \phi)} \sqrt{\mathcal{E}_1(U_1(1_{F_\ell}\mu_n), U_1(1_{F_\ell}\mu_n))}. \label{eq:4.8-2}
\end{equation}
If $\mathcal{E}_1(U_1(1_{F_\ell}\mu_n), U_1(1_{F_\ell}\mu_n)) = 0$ then the desired estimate is immediate. Otherwise, dividing both sides of \eqref{eq:4.8-1} by $\sqrt{\mathcal{E}_1(U_1(1_{F_\ell}\mu_n), U_1(1_{F_\ell}\mu_n))}$, we obtain
\begin{equation}
\sqrt{\mathcal{E}_1(U_1(1_{F_\ell}\mu_n), U_1(1_{F_\ell}\mu_n))} \le \sqrt{\mathcal{E}_1(U_1(1_{F_\ell}\mu), U_1(1_{F_\ell}\mu))} 
+ \varepsilon \sqrt{\mathcal{E}_1(\phi, \phi)}. \label{eq:4.8-3}
\end{equation}
Since $1_{F_\ell}\mu \in {\mathcal S}_0$, its energy $\mathcal{E}_1(U_1(1_{F_\ell}\mu), U_1(1_{F_\ell}\mu))$ is finite. By \eqref{eq:4.8-2} and \eqref{eq:4.8-3}, we have 
$$
\sup_n \mu_n(F_\ell) \le \sqrt{\mathcal{E}_1(\phi, \phi)} \left( \sqrt{\mathcal{E}_1(U_1(1_{F_\ell}\mu), U_1(1_{F_\ell}\mu))} 
+ \varepsilon \sqrt{\mathcal{E}_1(\phi, \phi)} \right)< \infty.
$$ 

We see that for each $\ell$, 
\begin{align*}
\lefteqn{\form_1(U_1(1_{F_\ell}\mu_n)-U_1(1_{F_\ell}\mu),U_1(1_{F_\ell}\mu_n)-U_1(1_{F_\ell}\mu))}\\
=& \ \int_{F_\ell} \Big(U_1(1_{F_\ell}\mu_n)(x)-U_1(1_{F_\ell}\mu)(x)\Big)\mu_n(dx)  -
\int_{F_\ell} \Big(U_1(1_{F_\ell}\mu_n)(x)-U_1(1_{F_\ell}\mu)(x)\Big)\mu(dx)  \\
=  & \ \int_{F_\ell} \Big({\mathbb E}_x\Big[\int_0^\infty   e^{-s}d(1_{F_\ell}{\sf A}^n)_s\Big]  -
{\mathbb E}_x\Big[\int_0^\infty e^{-s}d(1_{F_\ell}{\sf A})_s\Big]  \Big) \mu_n(dx) \\
& \qquad 
 - \int_{F_\ell} \Big({\mathbb E}_x\Big[\int_0^\infty  e^{-s}d(1_{F_\ell}{\sf A}^n)_s\Big]  -
{\mathbb E}_x\Big[\int_0^\infty  e^{-s}d(1_{F_\ell}{\sf A})_s\Big]  \Big) \mu(dx) \\
=:& \  {\sf (I)} - {\sf (II)}. 
\end{align*}

We proceed to estimate the term {\sf (I)}.\begin{align*}
|{\sf (I)}| & \le \Big| \int_{F_\ell} \Big({\mathbb E}_x\Big[\int_0^t  e^{-s}d(1_{F_\ell}{\sf A}^n)_s\Big]  -
{\mathbb E}_x\Big[\int_0^t e^{-s}d(1_{F_\ell}{\sf A})_s\Big]  \Big) \mu_n(dx)\Big|   \\
& \quad + \int_{F_\ell} {\mathbb E}_x\Big[\int_t^\infty  e^{-s}d(1_{F_\ell}{\sf A}^n)_s\Big]  \mu_n(dx) + 
\int_{F_\ell} {\mathbb E}_x\Big[\int_t^\infty e^{-s}d(1_{F_\ell}{\sf A})_s\Big] \mu_n(dx)   \\
& \le \sup_{x\in F_\ell}
\Big|
{\mathbb E}_x \Big[ \int_0^t e^{-s} d(1_{F_\ell}{\sf A}^n)_s\Big]-
{\mathbb E}_x \Big[ \int_0^t e^{-s} d(1_{F_\ell}{\sf A})_s\Big]\Big| \mu_n(F_\ell)  \\
& \quad + 2  \sup_{x\in F_\ell,\, n \in{\mathbb N}} 
\Big( {\mathbb E}_x \Big[ \int_t^\infty e^{-s} d(1_{F_\ell}{\sf A}^n)_s \Big]  \vee
 {\mathbb E}_x \Big[ \int_t^\infty e^{-s}d(1_{F_\ell}{\sf A})_s \Big] \Big) \mu_n(F_\ell) \\
& \le  \sup_n \mu_n(F_\ell) \, \sup_{x\in F_\ell}
\Big|
{\mathbb E}_x \Big[ \int_0^t e^{-s} d(1_{F_\ell}{\sf A}^n)_s\Big]-
{\mathbb E}_x \Big[ \int_0^t e^{-s} d(1_{F_\ell}{\sf A})_s\Big]\Big|   \\
& \quad + 2\,  \sup_n \mu_n(F_\ell) \,  \sup_{x\in F_\ell,\, n \in{\mathbb N}} 
\Big( {\mathbb E}_x \Big[ \int_t^\infty e^{-s} d(1_{F_\ell}{\sf A}^n)_s \Big]  \vee
 {\mathbb E}_x \Big[ \int_t^\infty e^{-s}d(1_{F_\ell}{\sf A})_s \Big] \Big) 
\end{align*}
and the right hand side goes to $0$ when $n\to \infty$ and then $t \to \infty$ by assumptions (1) and (2). 
The term {\sf (II)} can be estimated similarly to {\sf (I)}. Following the same steps, we find that it also goes to 
$0$ as $n\to\infty$. Thus it follows that 
$$
\lim_{n\to \infty} \form_1(U_1(1_{F_\ell}\mu_n)-U_1(1_{F_\ell}\mu))=0
$$
for each $\ell$, and whence
\begin{align*}
\rho(\mu_n, \mu) 
& =  \sum_{\ell=1}^\infty \frac 1{2^\ell} \Big(1\wedge \rho_0(1_{F_\ell}\mu_n, 1_{F_\ell}\mu)\Big) \\
& = \sum_{\ell=1}^\infty \frac 1{2^\ell} \Big(1\wedge \sqrt{\form_1\left(U_1(1_{F_\ell}\mu_n)-U_1(1_{F_\ell}\mu)\right)}\Big) \to 0 \quad (n\to \infty).
\end{align*}\end{proof}

\begin{exam} \label{ex-4.1}
We revisit the case of $(\form, \dom)=(\frac 12{\mathbb D}, H^1(\real^d))$ and $d\ge 2$.  
Take $\beta \in \real$, $a\in \real^d$, and sequences $\{\beta_n\} \subset \real$ and $\{a_n\} \subset \real^d$ 
such that $\beta_n \to \beta$ (in $\real$) and $a_n \to a$ (in $\real^d$).

Set $\mu_n(dx):=|x-a_n|^{-\beta_n}dx$ and $\mu(dx):=|x-a|^{-\beta}dx$. These measures are known to be smooth.  Since $a_n$ converges to $a$, we can choose a subsequence $\{n_\ell\}$ such that $n_\ell < n_{\ell+1}$ for each $\ell \in\mathbb{N}$ and $|a_k-a|<1/\ell$ for all $k\ge n_\ell$.  Using the sequence $\{a_n\}$ and the subsequence $\{n_\ell\}$, we define the nest $\{F_\ell\}_{\ell \in \mathbb{N}}$ as follows:
$$
F_\ell := \big\{ x\in\real^d : 1/\ell \le |x-a|\le \ell \big\} \setminus 
\Big(\bigcup_{k=1}^{n_\ell-1} \big\{x\in \real^d : |x-a_k|<1/\ell \big\} \Big). 
$$
Then, it is straightforward to verify that $\{\mu_n, \mu\}\subset\mathcal{S}(\{F_\ell\})$ and that the conditions 
of the Theorem \ref{conv_rho_nest} are satisfied for the nest $\{F_\ell\}$.  
\end{exam}

\section{Examples} \label{sec_example}

In this section, we present several examples to illustrate the results discussed thus far.

\begin{exam} \label{ex-03}
We consider a $d$-dimensional Brownian motion  ${\mathbb M}$ for
$d\ge 2$ and its Dirichlet form $(\mathcal{E}, \mathcal{F}) = (\frac{1}{2}{\mathbb D}, H^1(\mathbb{R}^d))$ on $L^2(\mathbb{R}^d;m)$, where $m(dx)$ denotes the Lebesgue measure. (see Examples \ref{ex-01} and \ref{ex-4.1}). Let $\mu$ be the Lebesgue measure on a $(d-1)$-dimensional hyperplane $F\subset \real^d$.  
 Since the Brownian motion on $\real^d$ is rotationally invariant, we may assume  without loss of generality that
 $F$ is the hyperplane $F=\{(x^{(1)}, 0) : x_{(1)} \in \mathbb{R}^{d-1}\}$.
We define a nest of sets relevant to this hyperplane as follows: for each $\ell \in{\mathbb N}$, 
let $F_\ell$ be the closed ball in $\real^d$ with radius $\ell$, centered at the origin, i.e., 
$F_{\ell} := \{ x\in \real^d: |x| \le \ell\}$.

 It is clear that $\mu(F_{\ell}) < \infty$ for any $\ell$ because  the intersection between 
 $F$ and $F_{\ell}$ is compact.   For any compact set $K \subset {\mathbb R}^d$, there exists
 a number $\ell_0$ such that $K \subset F_{\ell_0}$, so
 ${\rm Cap}(K \setminus F_{\ell}) = 0$ for any $\ell \ge \ell_0$.
 Hence we see that the measure $\mu$ is smooth.
 Using \cite[Exercise 2.2.1]{FOT11}, we find that
 $1_{F_{\ell}} \mu \in S_0$ for all $\ell$.
 Thus, $\mu \in S(\{F_{\ell}\})$.

 Let $\{E_n\}$ be any nest. We define $\mu_n = 1_{E_n} \mu$.
It is then clear that $\{\mu_n, \mu\} \subset S(\{F_{\ell}\})$.
 Let ${\sf A}^n, {\sf A}$ be PCAFs corresponding to
 $\mu_n, \mu$, respectively.

 We now verify assumption {\sf (A1)} from Theorem \ref{thm:main}.
For any compact set $K$, we can bound the quantity as follows:
 \begin{align*}
  &\hspace{-0.5cm} \sup_{0 \le s \le t} \int_{K}
  \left|{\mathbb E}_x\left[(1_{F_{\ell}} {\sf A}^n)_s\right] -
  {\mathbb E}_x\left[(1_{F_{\ell}} {\sf A})_s\right]\right| m(dx) \\
  &\le \sup_{0 \le s \le t} \int_{K}
  {\mathbb E}_x\left[|(1_{F_{\ell}} {\sf A}^n)_s
  - (1_{F_{\ell}} {\sf A})_s|\right] m(dx) \\
  &\le \sup_{0 \le s \le t} \int_{K}
  {\mathbb E}_x\left[\int_{0}^{s} \left|1_{F_{\ell}}(X_u)
  1_{E_n}(X_u) - 1_{F_{\ell}}(X_u)\right| d{\sf A}_u\right] m(dx) \\
  &\le \int_{K} {\mathbb E}_x
  \left[\int_{0}^{t} 1_{F_{\ell} \cap E_n^c} (X_s) d{\sf A}_s\right] m(dx)\\ 
  &\le \int_{K} {\mathbb E}_x\left[\int_{\sigma_{E_n^c}}^{t} 1_{E_n^c}
  (X_s) d{\sf A}_s\right] m(dx)\\
  & \le \int_{K} {\mathbb E}_x[{\sf A}_t;  \sigma_{E_n^c} < t] m(dx).
 \end{align*}
 Since \(d\mu(x)=dx^{(1)}d\delta_0(x^{(2)})\) for \(x=(x^{(1)}, x^{(2)})\in \mathbb{R}^{d-1}\times \mathbb{R}\) where \(\delta_0\) is a Dirac measure at \(0\), \({\sf A}\) has the same law as a local time at \(0\) for a one-dimensional Brownian motion \(W\) (see \cite[\S 5.3.\((4^\circ)\)]{CF12}). Therefore, by Tanaka's formula, for \(a:= |x^{(1)}|\), we have
\[\mathbb{E}_x[{\sf A}_t]=\mathbb{E}_{a}|W_t|-|a| \le \mathbb{E}_{a}|W_t-a| =\int_0^\infty \frac{e^{-\frac{|y-a|^2}{2t}}}{\sqrt{2\pi t}}|y-a| dy \leq 2\int_0^\infty \frac{e^{-\frac{|y|^2}{2t}}}{\sqrt{2\pi t}}y dy   = \sqrt{\frac{2t}{\pi}}.\]
 Since $\{E_n\}$ is a nest, we have  $\mathbb{P}_x(\sigma_{E_n^c} <t) \to 0$ as \(n\to \infty\). By the dominated convergence theorem, this implies 
$$\lim_{n\to \infty} \int_K \mathbb{E}_x[{\sf A}_t; \sigma_{E_n^c} < t]  m(dx)=0.$$
Thus, assumption {\sf (A1)} is satisfied.

  Next we check {\sf (A2)}.  Since \(\mu_n\) is increasing,  it is sufficient to prove that
  for any compact set $K \subset \bR^d$,
  \begin{align*}
   \lim_{t \to \infty} \int_{K} \bE_x\left[\int_{t}^{\infty}
   e^{-\alpha s} d(1_{F_{\ell}} {\sf A})_s\right] m(dx)= 0.
  \end{align*}
 Recalling that 
 $\mathbb{E}_x[{\sf A}_t] \le \sqrt{2t/\pi}$ and \(\zeta=\infty\), we can apply integration by parts 
 to the integral inside the expectation:
  \begin{align*}
   \int_{K} \bE_x\left[\int_{t}^{\infty} e^{-\alpha s}
   d(1_{F_{\ell}}{\sf A})_s\right]
   &= \int_{K} \bE_x\left[\left[e^{-\alpha s}
   (1_{F_{\ell}}{\sf A})_s\right]_t^{\infty}
   + \alpha\int_{t}^{\infty} e^{-\alpha s}
   (1_{F_{\ell}} {\sf A})_s\, ds\right] m(dx)\\
   &= \int_{K} \bE_x\left[-e^{-\alpha t} (1_{F_{\ell}}{\sf A})_t
   + \alpha \int_{t}^{\infty} e^{-\alpha s}
   (1_{F_{\ell}}{\sf A})_s\, ds\right] m(dx)\\
   &\le \alpha \int_{K} \int_{t}^{\infty}
   e^{-\alpha s} \mathbb{E}_x[{\sf A}_s]\, ds m(dx)\\
   &\le \alpha m(K) \int_{t}^{\infty}
   e^{-\alpha s} \sqrt{\frac{2s}{\pi}}\, ds
   \xrightarrow{t \to \infty} 0.
  \end{align*}
  Therefore, $\{{\sf A}^n, {\sf A}\}$ satisfies assumption (A2). Hence, by Theorem \ref{thm:main}, \(\mu_n\) converges to \(\mu\) vaguely on the nest \(\{F_{\ell}\}\) in the resolvent sense. This means that we approximated 
a singular Radon measure not of finite energy integrals by a  sequence of singular Radon measures of finite energy integrals.
\qed
\end{exam}

\begin{exam}[{\small \bf Approximation of a Nowhere Radon Measure by Singular Spherical Measures}]\label{ex-02} \ 
In this example, we introduce {\it a nowhere Radon smooth measure}. Such a measure, 
whose restriction to any non-empty open set has infinite volume, was constructed by 
Stollmann-Voigt \cite{SV85} and later by Albeverio-Ma \cite{AM92}.  

We consider again a $d$-dimensional Brownian motion 
$(d\ge 2$) and its Dirichlet form $(\mathcal{E}, \mathcal{F}) = (\frac{1}{2}{\mathbb D}, H^1(\mathbb{R}^d))$ on $L^2(\mathbb{R}^d;m)$, where $m(dx)$ denotes the Lebesgue measure. Our goal is to realize a nowhere Radon measure as the limit of a sequence of singular measures, each supported on the surfaces of countably many concentric spheres. Although this sequence oscillates significantly and lacks monotonicity, which makes its vague convergence analytically intricate, the probabilistic behaviour of the corresponding PCAFs enables a seamless application of Theorem \ref{thm:main}.

In what follows, we outline a construction of a nowhere Radon measure and its associated compact nest, and then consider the convergence of measures. Let $\{x_j\}$ be a dense subset of ${\mathbb R}^d$ and  $E_{\ell} = \{x \in {\mathbb R}^d : |x| \le \ell\}$. 
 Since each point $x_j$ has zero capacity, we can find a sequence of open balls $\{G_{j,k}\}$ centered at
 $x_j$ that is decreasing in $k$ and satisfies $\bigcap_{k \ge 1} G_{j,k} = \{x_j\}$ and 
 ${\rm Cap}(G_{j,k}) + m(G_{j,k})  \le 2^{-k}$. We take $\alpha_j > d$ for each $j$. Let the sets $\{G_\ell\}$ and $\{F_\ell\}$ be defined as follows:
 \begin{align*}
  G_{\ell} = \left(\bigcup_{j=1}^{\ell} G_{j,2\ell}\right)
  \cup \left(\bigcup_{j=\ell+1}^{\infty}  G_{j,j}\right) \quad
  \text{and} \quad F_{\ell} = E_{\ell} \setminus G_{\ell}
 \end{align*}
and we choose a sequence of positive numbers $\{c_j\}$ that satisfies 
$\dis c_j \sup_{E_j \setminus G_{j,j}} |x - x_j|^{-\alpha_j} \le 2^{-j}$.

Now, we define a function $f$ and a measure $\mu$. Let $f$ be given by
$\dis f(x) = \sum_{j \ge 1} c_j |x - x_j|^{-\alpha_j}$. 
The measure $\mu$ is then defined as 
$\dis \mu(dx) = f(x) m(dx).$ By \cite[Proposition 1.3]{AM92}, $\mu=f \cdot m$ is a smooth measure on ${\mathbb R}^d$ whose generalized nest is $\{L_{\ell}\}_{\ell \ge 1}$. We note that  \(\mu(O)= \infty\) holds for any non-empty open set \(O\) because $\alpha_j > d$.

Since $F_{\ell}$ is compact and $f(x)$ is
  bounded continuous on $F_{\ell}$ for any $\ell \in \mathbb{N}$,
  we see that for any $v \in \cF \cap C_{0}(\real^d)$,
  \begin{align}
  \label{eq:ex_T_1}
    \int_{\real^d} |v(x)| 1_{F_{\ell}}\mu(dx) &= \int_{\real^d}  |v(x)| 1_{F_{\ell}} f(x) m(dx) 
    \le \|1_{F_{\ell}}f\|_{L^{2}} \|v\|_{L^{2}} \le \|1_{F_{\ell}}f\|_{L^{2}} \sqrt{\cE_{1}(v,v)}. 
  \end{align}
  Hence, the nest $\{F_{\ell}\}$ is also an attached nest of $\mu$,
  that is, $\mu \in \mathcal{S}(\{F_{\ell}\})$.

Next let $\sigma_{x_j, r}$ denote the uniform surface measure on the sphere $\{x \in \mathbb{R}^d : |x - x_j| = r\}$, normalized such that $\int_0^\infty \sigma_{x_j, r}(dx) dr = m(dx)$. Since $\sigma_{x_j, r}$ are singular, we define the sequence of singular measures $\mu_n$ by
\begin{align*}
    \mu_n(dx) := \frac{1}{n} \sum_{j=1}^n c_j \sum_{k=1}^{n^2} \left( \frac{k}{n} \right)^{-\alpha_j} \sigma_{x_j, \frac{k}{n}}(dx).
\end{align*}
By construction, each $\mu_n$ is supported exclusively on a collection of thin spherical surfaces. Since \(\sigma_{x_j, \frac{k}{n}}\) are Radon measures of finite energy integrals by the same argument  as \cite[Example 5.3.\((3^\circ)\)]{CF12}, so is \(\mu_n\). As $n$ increases, the radii $r_{k,n} = k/n$ shift entirely, rendering the sequence $\mu_n$ heavily oscillating. Since the sequence is neither monotonic nor absolutely continuous, verifying the vague convergence of $\{\mu_n\}$ to $\mu$ through direct analytical methods is quite formidable.

However, our framework provides a robust alternative to resolve this difficulty. Let ${\sf A}^n_t$ and ${\sf A}_t$ be the PCAFs associated with $\mu_n$ and $\mu$, respectively. The PCAF ${\sf A}^n_t$ is a weighted sum of the \textit{spherical local times} of the Brownian motion on the concentric spheres. 

To apply Theorem \ref{thm:main}, we analyze the behaviour of the restricted measures $1_{F_\ell}\mu_n$ on the nest and check the conditions {\sf (A1)} and {\sf (A2)} in Theorem \ref{thm:main}.

First we check {\sf (A2)} by showing \eqref{1-pot}. Recall that $F_\ell = E_\ell \setminus G_\ell$.  
 If a point $y$ belongs to $F_\ell$, it must lie outside $G_\ell$. Let $r_{j,k}$ be the radius of $G_{j,k}$. 
Then we have \(|y - x_j| \ge r_{j, 2 \ell}\) for \( j \le \ell\) and \(|y - x_j| \ge r_{j, j}\) for \(j \ge \ell+1.\) 
Let \(R_1(x,y)\) be a \(1\)-order resolvent density kernel for Brownian motion, 
and we set \(M_{\ell, j} := c_j (r_{j, 2 \ell})^{-\alpha_j}\) for $j \le \ell$. Since $c_j (r_{j, j})^{-\alpha_j} \le 2^{-j}$ for \(j\ge \ell+1\), 
we have the following bound for the \(1\)-order potential \(U_1(1_{F_\ell}\mu_n)\) of \(1_{F_\ell}\mu_n\).
\begin{align}
& U_1(1_{F_\ell}\mu_n)(x) = \int_{F_\ell} R_1(x, y) \mu_n(dy)=\sum_{j=1}^n \sum_{k=1}^{n^2} \frac{1}{n} \int_{F_\ell} R_1(x, y) \frac{c_j}{|y-x_j|^{\alpha_j} }\sigma_{x_j, \frac{k}{n}}(dy) \\
&\leq \sum_{j=1}^\ell M_{\ell, j} \sum_{k=1}^{n^2} \frac{1}{n} \int_{F_\ell} R_1(x, y) \sigma_{x_j, \frac{k}{n}}(dy) + \sum_{j=\ell+1}^n 2^{-j} \sum_{k=1}^{n^2} \frac{1}{n} \int_{F_\ell} R_1(x, y) \sigma_{x_j, \frac{k}{n}}(dy) . \label{eq:u1_bound_sum}
\end{align}
Let $H_x(r) := \int_{\mathbb{R}^d} R_1(x, y) \sigma_{x_j, r}(dy)$ for $r > 0$. By the polar coordinate decomposition of the Lebesgue measure $m(dy)$, its integral is bounded: $\int_0^\infty H_x(r) dr = R_1 1(x) \le 1$. Since the Riemann sum $\sum_{k=1}^{n^2} \frac{1}{n} H_{x_j}(\frac{k}{n})$ converges to the integral $\int_0^\infty H_x(r) dr$, by \eqref{eq:u1_bound_sum}, we obtain
\begin{align*}
    U_1(1_{F_\ell}\mu_n)(x) \le C \left( \sum_{j=1}^\ell M_{\ell, j} + \sum_{j=\ell+1}^n 2^{-j} \right) 
    \le C \left( \sum_{j=1}^\ell M_{\ell, j} + 2^{-\ell} \right) =: \tilde{M}_\ell < \infty.
\end{align*}
for some \(C>0\).
This bound $\tilde{M}_\ell$ is independent of $n$ and $x \in F_\ell$ for each \(\ell\), and so \eqref{1-pot} holds. Thus, Proposition \ref{propforA2} guarantees that the tail estimate condition {\sf (A2)} automatically holds.

Next we check {\sf (A1)}. To this end, fix any compact set $K \subset \mathbb{R}^d$, $\ell \in{\mathbb N}$ and $t > 0$.  
We first verify the pointwise convergence of the expectations. The limit measure $\mu(dy) = f(y)m(dy)$ is defined by the density 
$f(y) = \sum_{j \ge 1} c_j |y - x_j|^{-\alpha_j}$.  While $\mu$ has this explicit density, the approximating measures $\mu_n$ do not 
possess density functions, which makes a direct evaluation of the difference of PCAFs intricate.  However, drawing on the convergence 
of the Riemann sum demonstrated in the proof of {\sf (A2)} and using the transition density $p_u(x, y)$ of the Brownian motion, 
we can translate the expectation of the PCAF into an integral with respect to the spatial measure.  For each $s \in (0, t]$ and $x \in K$, we set 
\begin{align*}
\Psi_j(r, s) := \int_{F_\ell} \left(\int_0^s p_u(x, y)du\right) \frac{c_j}{|y-x_j|^{\alpha_j}} \sigma_{x_j,r}(dy).
\end{align*}
To rigorously justify the interchange of the limit $n \to \infty$ and the infinite sum over $j$, we fix an arbitrary integer $J > \ell$ 
and evaluate the pointwise difference:
\begin{align} \label{eq:6.2}
&\!\!\! \left|{\mathbb E}_x[(1_{F_\ell}{\sf A}^n)_s]- {\mathbb E}_x[(1_{F_\ell}{\sf A})_s] \right| \notag \\
&= \left| \sum_{j=1}^n \frac{1}{n} \sum_{k=1}^{n^2} \Psi_j\left(\frac{k}{n}, s\right) - \sum_{j=1}^\infty \int_0^\infty \Psi_j(r, s)dr \right| \notag \\
&\le \sum_{j=1}^J \left| \frac{1}{n} \sum_{k=1}^{n^2} \Psi_j\left(\frac{k}{n}, s\right) - \int_0^n \Psi_j(r, s) dr \right| 
+ \sum_{j=J+1}^n \frac{1}{n} \sum_{k=1}^{n^2} \Psi_j\left(\frac{k}{n}, s\right)  \\
&\quad + \sum_{j=1}^J \int_n^\infty \Psi_j(r, s)dr + \sum_{j=J+1}^\infty \int_0^\infty \Psi_j(r, s)dr. 
\end{align}
We evaluate the four terms on the right-hand side of \eqref{eq:6.2} separately.  

For the first term, similar to $H_x(r)$ in the proof of {\sf (A2)}, for each fixed $j \le J$, 
the mapping $r \mapsto \Psi_j(r, s)$ is compactly supported in $(0, \infty)$ and Riemann integrable since the boundary of 
$F_\ell$ consists of spheres. Thus, the discrete sum $\frac{1}{n} \sum_{k=1}^{n^2} \Psi_j(k/n, s)$ is exactly the finite 
Riemann sum approximation of the continuous integral. Consequently, the first term converges to 0 as $n \to \infty$.

We apply the bound $c_j |y-x_j|^{-\alpha_j} \le 2^{-j}$ for $j \ge \ell+1$ on the set $F_\ell$ to estimate the second term. 
Following the same logic as in {\sf (A2)} 
with the heat kernel integral $\int_0^s p_u(x, y)du$ in place of the resolvent $R_1(x, y)$, its corresponding Riemann sum is uniformly bounded 
by some $C' s$.  Hence, by extending the finite sum up to $n$ to an infinite sum, the second term is bounded by 
$C' s \sum_{j=J+1}^\infty 2^{-j} = C' s 2^{-J}$.

For the third term, since the integrals over the whole space are finite, i.e., 
$\int_0^\infty \Psi_j(r, s)dr \le e^s U_1(1_{F_\ell}\mu)(x) < \infty$, the tails of these finite $J$ integrals converge to 0 as $n \to \infty$.

For the fourth term, similar to the second term, its continuous integral with respect to $\sigma_{x_j, r}(dy) dr$ is bounded 
by $\int_0^s \int_{\mathbb{R}^d} p_u(x,y)dydu \le s$. Thus, the fourth term is bounded by $s \sum_{j=J+1}^\infty 2^{-j} = s 2^{-J}$. 

Notice that the sum of the upper bounds for the second and fourth terms is $(C'+1) s 2^{-J}$, which is independent of $n$ and can be made 
arbitrarily small by choosing $J$ sufficiently large. Therefore, taking the limit superior as $n \to \infty$ for the pointwise difference yields:
\begin{align*}
\limsup_{n \to \infty} \left|{\mathbb E}_x[(1_{F_\ell}{\sf A}^n)_s]- {\mathbb E}_x[(1_{F_\ell}{\sf A})_s] \right| \le (C'+1)s 2^{-J}. 
\end{align*}
Since $J \ge \ell+1$ can be chosen arbitrarily large, taking the limit $J \to \infty$ shows that the right-hand side vanishes.  
This establishes the pointwise convergence $\lim_{n \to \infty} {\mathbb E}_x[(1_{F_\ell}{\sf A}^n)_s] = {\mathbb E}_x[(1_{F_\ell}{\sf A})_s]$ 
for each $s \in [0, t]$.

Furthermore, for each fixed $x \in K$, the mapping $s \mapsto f_n(s):={\mathbb E}_x[(1_{F_\ell}{\sf A}^n)_s]$ is  monotonically increasing 
with respect to $s$ for each $n$, and the limit function $s \mapsto f(s):={\mathbb E}_x[(1_{F_\ell}{\sf A})_s]$ is also monotonically increasing 
and continuous on $[0, t]$.  Then we find that, for each $x\in K$,
\begin{align*}
\lim_{n\to\infty} \sup_{0\le s\le t} \left| {\mathbb E}_x[(1_{F_\ell}{\sf A}^n)_s] - {\mathbb E}_x[(1_{F_\ell}{\sf A})_s] \right| = 0.
\end{align*}
Finally, integrating this uniform difference over the compact set $K$ with respect to $m(dx)$, and combining it with the uniform bound $\sup_{0\le s\le t} {\mathbb E}_x[(1_{F_\ell}{\sf A}^n)_s] \le e^t U_1(1_{F_\ell}\mu_n)(x) \le e^t \tilde{M}_\ell$, we can apply the dominated convergence theorem to conclude that
\begin{align*}
&\!\!\! \lim_{n\to\infty} \sup_{0\le s\le t} \int_K \left| {\mathbb E}_x[(1_{F_\ell}{\sf A}^n)_s] - {\mathbb E}_x[(1_{F_\ell}{\sf A})_s] \right| m(dx) \\
&\le \lim_{n\to\infty} \int_K \sup_{0\le s\le t} \left| {\mathbb E}_x[(1_{F_\ell}{\sf A}^n)_s] - {\mathbb E}_x[(1_{F_\ell}{\sf A})_s] \right| m(dx) \\
&= 0,
\end{align*}
which implies that {\sf (A1)} holds. \hfill \fbox{}
\end{exam}

\appendix
\section{Characterization and Continuity of the Revuz correspondence through the Killing Procedure}
\label{sec-killing}
Although not directly related to the convergence results in Section 4, this appendix provides an alternative perspective on the continuity of the Revuz correspondence.
While Theorem \ref{thm-Revuz} characterizes the Revuz correspondence through the underlying measure $m$, 
this section presents a different perspective.  We provide a new characterization by analyzing processes with 
killing inside the domain or continuous killing at the boundary with the boundary itself  identified 
with the cemetery point $\partial$.

Recall that  ${\mathbb M}=(X_t, {\mathbb P}_x)$ is an $m$-symmetric Hunt process whose Dirichlet form $(\form, \dom)$ is assumed to be regular and irreducible.

\begin{thm}[{The Beurling--Deny decomposition, c.f. \cite[Theorem 3.2.1]{FOT11}, \cite[Section 4.3]{CF12}}] \label{BDdecomp}
 There exist a strongly local Dirichlet form \(\mathcal{E}^{(c)}\), a symmetric Radon measure \(J\) on 
 \((E\times E) \setminus {\sf diag}\), called the jumping measure,  where ${\sf diag}=\{(x,x) : x\in E\}$ 
 denotes the diagonal set, and a Radon measure  $\kappa$  on \(E\),  called the killing measure,  
 such that, for \(f,g\in \mathcal{F}\), 
\[\mathcal{E}(f,g)=\mathcal{E}^{(c)}(f,g)+\frac{1}{2}\iint_{E\times E \setminus {\sf diag}} (f(x)-f(y)) (g(x)-g(y))J(dxdy)+\int_E fg d\kappa.\]
\end{thm}

The killing measure $\kappa$ from the Beurling--Deny decomposition is a smooth measure.  By considering the potentials of the killing measure restricted to an attached nest, we can derive the following equations, which are analogous to those presented in Theorem \ref{thm-Revuz}. 
We set \(\mathbb{E}_{\kappa}[\cdot]:= \int_E \mathbb{E}_x[\cdot] \kappa(dx).\)
For \(\alpha \geq 0\), we define \(\varphi_{\alpha}(x):=\mathbb{E}_x[e^{-\alpha\zeta} {1}_{\{\partial\}}(X_{\zeta -})]\) for \(x \in E\) and \(\varphi_{\alpha}(\partial):=0\), and we set \(\varphi:=\varphi_{0}\).

\begin{lem}
For \(\alpha, \beta \geq 0\), it holds that
\begin{equation}
\varphi_\alpha(x)-\varphi_\beta(x)+(\alpha-\beta)R_\beta\varphi_\alpha(x) = 0. \label{eq:resolvent_nu0}
\end{equation}
\end{lem}

\begin{proof}
We may assume \(\alpha >0\) or \(\beta>0\). Since \(1_{\{\partial\}}(X_{\zeta -}) \circ \theta_t = 1_{\{\partial\}}(X_{\zeta -})\) for \(t<\zeta\), we have
\begin{align*}
\frac{\varphi_\alpha(x)-\varphi_\beta(x)}{\beta -\alpha} &= \mathbb{E}_x\left[\int_0 ^\zeta e^{-\alpha \zeta} e^{-(\beta -\alpha)t}1_{\{\partial\}}(X_{\zeta -}) \,dt \right]\\
&= \int_0 ^\infty e^{-\beta t} \mathbb{E}_x\left[\left(e^{-\alpha \zeta}1_{\{\partial\}}(X_{\zeta -})\right) \circ \theta_t \, 1_{\{t<\zeta\}} \right]\,dt \\
&= \int_0 ^\infty e^{-\beta t}\mathbb{E}_x\left[\varphi_\alpha(X_t)\, 1_{\{t<\zeta\}}\right]\,dt \ =\ R_\beta \varphi_\alpha(x).
\end{align*}
\end{proof}

\begin{lem}\label{Revuzkappa}
For \(\mu \in \mathcal{S}\), its corresponding  PCAF \(\sfA\),  any \(t, \alpha \geq 0\) and 
\(f \in \mathcal{B}_+(E)\), it holds that
\begin{equation}\label{Revuz-kappa1}
\mathbb{E}_{\kappa} \left[ \int_0^\infty e^{-\alpha s} f(X_s) \, d\sfA_s  \right] 
= \int_E ( 1_E - \alpha R_\alpha  1_E-\varphi_{\alpha})f \, d\mu 
\end{equation}
and
\begin{equation}\label{Revuz-kappa2}
\mathbb{E}_{\kappa} \left[ \int_0^t f(X_s) \, d\sfA_s \right] 
= \int_E ( 1_E-\varphi-p_t( 1_E-\varphi))f\, d\mu.
\end{equation}
\end{lem}

\begin{proof}
First, we assume that \(\alpha>0\) and \(\mu \in \mathcal{S}_0\). Since \(\kappa \in \mathcal{S}\), there exists a nest $\{F_\ell\}$ such that $1_{F_\ell} \kappa \in \mathcal{S}_{0}$, \(U_{\alpha}(1_{F_\ell}\kappa)\) is bounded for each \(\ell\), and \(\kappa (E\setminus (\cup_\ell F_\ell))=0\). By using the L\'{e}vy system, we have \(U_{\alpha}({\bf 1}_{F_\ell} \kappa)(x) = \mathbb{E}_x[e^{-\alpha \zeta} {\bf 1}_{F_\ell}(X_{\zeta -})]\). See \cite[Theorem 4.2.1]{CF12} or the proof of \cite[Proposition 1.2]{Oo26} for details.  By \cite[Lemma 5.1.3]{FOT11}, for any bounded non-negative Borel measurable function $f$,
\(
\mathbb{E}_x \left[ \int_0^\infty e^{-\alpha s} f(X_s) \, d\sfA_s \right] 
\) is a quasi-continuous version of $U_\alpha(f \mu)$.  Hence, by the monotonicity and \cite[(3.5)]{Y97}, we have
\begin{align*}
& \mathbb{E}_{\kappa} \left[ \int_0^\infty e^{-\alpha s} f(X_s) \, d\sfA_s \right] = \lim_{\ell \to \infty} \int_E \mathbb{E}_x \left[ \int_0^\infty e^{-\alpha s} f(X_s) \, d\sfA_s \right]  {\bf 1}_{F_\ell}(x) d\kappa(x)\\
&= \lim_{\ell \to \infty}\int_E  U_{\alpha}(1_{F_\ell} \kappa)f d\mu = \int_E ( 1_E - \alpha R_\alpha  1_E-\varphi_{\alpha})f \, d\mu.
\end{align*}
Moreover, \eqref{Revuz-kappa1} with \(\alpha>0\) holds for any \(f \in \mathcal{B}_+(E)\) by approximating by bounded functions \(f\wedge n\). For any \(\mu \in \mathcal{S}\), we take a nest \(\{G_\ell\}\) such that \(1_{G_\ell} \mu \in \mathcal{S}_0\) and \(\mu(E\setminus (\cup _\ell G_\ell))=0.\) By \cite[Theorem 3.1.4 (ii)]{CF12}, \(\mathbb{P}_x(\lim_\ell \sigma_{G_\ell^c}<\zeta)=0\) for q.e. \(x\in E\), where \(\sigma_{B}\) is the first hitting time to a Borel set \(B.\) By the monotone convergence theorem, as \(\ell \to \infty\) in \eqref{Revuz-kappa1} for \(1_{G_\ell}\mu \), \eqref{Revuz-kappa1} holds for \(\mu \in \mathcal{S}\) and \(\alpha>0\).

For \(\alpha=0\), when the process \(X\) is not transient, by the irreducibility, \(X\) is recurrent and conservative, so \(\kappa = 0.\) Hence we may assume that \(X\) is transient. We take a finite measure \(\mu\in \mathcal{S}_0^{(0)}\) and a bounded positive Borel measurable function \(f\). Then \(f\mu \in \mathcal{S}_0^{(0)}\). By \cite[Exercise 2.2.4]{FOT11} (see also the proof of \cite[Theorem 2.4.2.(ii)]{FOT11}), we can take a nest $\{F_\ell\}$ satisfying that $1_{F_\ell} \kappa \in \mathcal{S}_{0}^{(0)}$, \(U(1_{F_\ell}\kappa)\) is bounded for each \(\ell\), and \(\kappa (E\setminus (\cup_\ell F_\ell))=0\). Then we have
\begin{align*}
& \mathbb{E}_{\kappa} \left[ \int_0^\infty f(X_s) \, d\sfA_s \right] = \lim_{\ell \to \infty} \int_E \mathbb{E}_x \left[ \int_0^\infty  f(X_s) \, d\sfA_s \right]  {\bf 1}_{F_\ell}(x) d\kappa(x)\\
&=  \lim_{\ell \to \infty}\mathcal{E}(U(1_{F_\ell} \kappa), U(f\mu)) = \lim_{\ell \to \infty}\int_E  U(1_{F_\ell} \kappa)f d\mu = \int_E (1_E - \varphi )f \, d\mu.
\end{align*}

Furthermore \eqref{Revuz-kappa1} with \(\alpha=0\) holds for any \(f \in \mathcal{B}_+(E)\) by approximating by bounded functions \(f\wedge n\). For any \(\mu \in \mathcal{S}\), by \cite[Exercise 2.2.4]{FOT11}, we take a nest \(\{G_\ell\}\) such that \(1_{G_\ell} \mu \in \mathcal{S}_0^{(0)}\) and \(\mu(E\setminus (\cup _\ell G_\ell))=0.\) By \cite[Theorem 3.1.4 (ii)]{CF12}, \(\mathbb{P}_x(\lim_\ell \sigma_{G_\ell^c}<\zeta)=0\) for q.e. \(x\in E\).  By the monotone convergence theorem, as \(\ell \to \infty\) in \eqref{Revuz-kappa1} for \(1_{G_\ell}\mu \), \eqref{Revuz-kappa1} holds for \(\mu \in \mathcal{S}\) and \(\alpha=0\).

Applying the Laplace transform, we obtain
\begin{align*}
& \int_0^\infty e^{-\alpha t} \mathbb{E}_{\kappa} \left[ \int_0^t f(X_s) \, d\sfA_s \right] dt = \mathbb{E}_{\kappa}\left[\int_0^\infty  \int_0^{\infty} {1}_{\{s\leq t\}} e^{-\alpha t} f(X_s) \, d\sfA_s dt\right]\\
&= \mathbb{E}_{\kappa}\left[\int_0^\infty  \int_s^{\infty} e^{-\alpha t} f(X_s) \, dt d\sfA_s\right] = \frac{1}{\alpha} \mathbb{E}_{\kappa} \left[ \int_0^\infty e^{-\alpha s}   f(X_s) \, d\sfA_s \right]
\end{align*}
and, by the identity \eqref{eq:resolvent_nu0}, we have
\begin{align*}
\lefteqn{\int_0^{\infty} e^{-\alpha t} \int_E f(x) ( 1_E-\varphi(x) -p_t( 1_E-\varphi)(x))\, d\mu(x) dt}\\
&= \frac{1}{\alpha} \int_E f(x)( 1_E-\alpha R_{\alpha}1_E -\varphi+\alpha R_{\alpha}\varphi)(x) \,d\mu(x)\\
&= \frac{1}{\alpha} \int_E f(x)( 1_E-\alpha R_{\alpha} 1_E -\varphi_{\alpha})(x) \,d\mu(x).
\end{align*}
Hence \eqref{Revuz-kappa2} follows from the uniqueness of Laplace transforms.
\end{proof}

\begin{rem}
 Assume that \(\kappa =hdm\). Then \(X\) is a killed process of a resurrected process \(X^{\rm res}\) by \(h dm\) (\cite[Theorem 5.2.17]{CF12}). We set  \({\sf A}_t^{{\rm res}} := \int_0^t h(X_s^{\rm res})\,ds\).
Since \(h(X_t^{{\rm res}})=h(\partial) = 0\) for \(t\geq \zeta^{{\rm res}}\) and \(1_{\{\partial\}}(X_{\zeta^{{\rm res}}-}^{{\rm res}})1_{\{\zeta^{{\rm res}} <\infty\}}=1_{\{\zeta^{{\rm res}} <\infty\}}\), we have
 \begin{align*}
  R_{\alpha}h(x) &= \mathbb{E}^{{\rm res}}_x\left[\int_0^{\infty} e^{-\alpha t}
				 e^{-{\sf A}_t^{{\rm res}} }
				 h(X_t^{\rm res})\,dt \right]\\
 &= \mathbb{E}^{{\rm res}}_x\left[\int_0^{\zeta^{{\rm res}}} e^{-\alpha t}
		  \frac{d}{dt}\left(-\exp{\left(-\int_0^t h(X_s^{\rm res})\,ds\right)} \right) dt \right]\\
 &= \mathbb{E}^{{\rm res}}_x \left[\left[- e^{-\alpha t}
e^{- {\sf A}_t^{{\rm res}}} \right]_{t=0}^{t=\zeta^{{\rm res}}} \right] - \mathbb{E}^{{\rm res}}_x\left[\int_0^{\zeta^{{\rm res}}}\alpha e^{-\alpha t}
e^{-{\sf A}_t^{{\rm res}}} dt \right]\\
&= 1_E(x)-\mathbb{E}^{{\rm res}}_x\left[ e^{-\alpha \zeta^{{\rm res}}}
e^{-{\sf A}_{\zeta^{{\rm res}}}^{{\rm res}}}  \right]-\mathbb{E}^{{\rm res}}_x\left[\int_0^\infty \alpha e^{-\alpha t}e^{-{\sf A}_t^{{\rm res}}}  1_E(X_t^{\rm res})dt \right]\\
&= 1_E(x)-\mathbb{E}^{{\rm res}}_x\left[ e^{-\alpha \zeta^{{\rm res}}}
e^{-{\sf A}_{\zeta^{{\rm res}}}^{{\rm res}}} 1_{\{\partial\}}(X_{\zeta^{{\rm res}}-}^{{\rm res}})1_{\{\zeta^{{\rm res}} <\infty\}} \right]-\alpha R_{\alpha}1_E(x)\\
&= 1_E(x)-\mathbb{E}_x\left[ e^{-\alpha \zeta} 1_{\{\partial\}}(X_{\zeta-}) 1_{\{\zeta  <\infty\}}\right]-\alpha R_{\alpha}1_E(x)\\
&= 1_E(x)-\varphi_\alpha(x)-\alpha R_{\alpha}1_E(x),
 \end{align*}
so \eqref{Revuz-kappa1} for \(\kappa =h dm\) does not contradict \eqref{eq-Revuz}. Here, 
the second to last equation is justified by using \cite[\S A.3.2]{CF12} as follows. Let \(\hat{\Omega} := \Omega \times [0,\infty]\) and write \(\hat{\omega}:=(\omega, \lambda)\). We define 
\[\hat{X}_t(\hat{\omega}) := \begin{cases}X_t^{{\rm res}}(\omega) &: t<\lambda, \\ \partial &: t \geq \lambda, \end{cases} \text{ \ and\ }\hat{\zeta}(\hat{\omega}):=\zeta^{{\rm res}}(\omega) \wedge \lambda. \]
For each \(\omega \in \Omega,\) we take a probability measure \(P^{\omega}\) on \(([0,\infty], \mathcal{B}[0,\infty))\) such that \(P^\omega(\lambda \geq t) = e^{-{\sf A}_t^{{\rm res}}(\omega)}.\) By \cite[Lemma A.3.12]{CF12}, \(\hat{X}\) has the same law as that of \(X\). Since \(X_{\lambda-}^{{\rm res}}(\omega)\in E\),
\begin{eqnarray*}
&&\mathbb{E}_x\left[ e^{-\alpha \zeta} 1_{\{\partial\}}(X_{\zeta-}) 1_{\{\zeta  <\infty\}}\right]= \mathbb{E}^{{\rm res}}_xE^{\omega}\left[ e^{-\alpha \zeta^{{\rm res}}(\omega)}
1_{\{\partial\}}(X_{\zeta^{{\rm res}}-}^{{\rm res}}(\omega))1_{\{\zeta^{{\rm res}}(\omega) <\infty\}}1_{\{\zeta^{{\rm res}}(\omega) <\lambda\}} \right]\\
&&\hspace{5mm}+\mathbb{E}^{{\rm res}}_xE^{\omega}\left[ e^{-\alpha\lambda}
1_{\{\partial\}}(X_{\lambda-}^{{\rm res}}(\omega))1_{\{\lambda <\infty\}}1_{\{\zeta^{{\rm res}}(\omega) \geq \lambda\}} \right]\\
&&= \mathbb{E}^{{\rm res}}_x\left[ e^{-\alpha \zeta^{{\rm res}}}
1_{\{\partial\}}(X_{\zeta^{{\rm res}}-}^{{\rm res}})1_{\{\zeta^{{\rm res}} <\infty\}}e^{-{\sf A}_{\zeta^{{\rm res}}}^{{\rm res}}}  \right].
\end{eqnarray*}
\end{rem}

The killing measure \(\kappa\) characterizes the killing inside of the process. On the other hand, by extending an energy functional with \(\varphi_\alpha(x):=\mathbb{E}_{x}[e^{-\alpha \zeta}1_{\{\partial\}}(X_{\zeta -})]\), the following functional \(\nu_0\) is introduced in \cite{Oo26} to describe the energy of the part of the process that continuously escapes to the cemetery point.

\begin{definition}\label{def_nu_0}
For any positive universally measurable function \(f,\) we define an energy functional \(\nu_0\) of \(f\) by
\[\int f d\nu_0 := \varlimsup_{s\searrow 0}\frac{1}{s} \int_E f(x)(\varphi_{\alpha}(x)-e^{-\alpha s}p_s\varphi_{\alpha}(x)) dm(x),\]
where \(\varphi_{\alpha}(x):=\mathbb{E}_x[e^{-\alpha \zeta}1_{\{\partial\}}(X_{\zeta -})]\) for \(\alpha \geq 0\).
\end{definition}

We prove that Definition \ref{def_nu_0} is well-defined. For convenience, we set \(\int f d\nu_0^{\alpha} := \varlimsup_{s\searrow 0}\frac{1}{s} \int_E f(x)(\varphi_{\alpha}(x)-e^{-\alpha s}p_s\varphi_{\alpha}(x)) dm(x)\). We take a positive Borel measurable function \(f\) and \(\alpha > \beta \geq 0\). Since \(\varphi_{\alpha}-e^{-\alpha s}p_s \varphi_{\alpha} = \mathbb{E}_x[e^{-\alpha \zeta}1_{\{\partial\}}(X_{\zeta -})1_{\{s\geq \zeta\}}]\), it holds that \(\int f d\nu_0^\alpha \leq \int f d\nu_0^\beta.\) Suppose that there exists \(s>0\) such that \(\int_E f \mathbb{E}_x[e^{-\beta \zeta}1_{\{\partial\}}(X_{\zeta -})1_{\{s\geq \zeta\}}] dm <\infty\).  By the subadditivity of the limit superior and the dominated convergence theorem, we have
\begin{eqnarray*}
\lefteqn{\int f d\nu_0^\beta  = \varlimsup_{s\searrow 0}\frac{1}{s} \int_E f \mathbb{E}_x[e^{-\beta \zeta}1_{\{\partial\}}(X_{\zeta -})1_{\{s\geq \zeta\}}] dm}\\
&\leq & \int f d\nu_0^\alpha  + \varlimsup_{s\searrow 0} \frac{1}{s} \int_E f \mathbb{E}_x[(e^{-\beta \zeta}-e^{-\alpha \zeta})1_{\{\partial\}}(X_{\zeta -})1_{\{s\geq \zeta\}}] dm \\
&\leq & \int f d\nu_0^\alpha  + \varlimsup_{s\searrow 0} \frac{1}{s} \int_E f \mathbb{E}_x[(\alpha-\beta)\zeta e^{-\beta \zeta} 1_{\{\partial\}}(X_{\zeta -})1_{\{s\geq \zeta\}}] dm \\
&\leq & \int f d\nu_0^\alpha  + \varlimsup_{s\searrow 0} (\alpha-\beta) \int_E f \mathbb{E}_x[e^{-\beta \zeta}1_{\{\partial\}}(X_{\zeta -})1_{\{s\geq \zeta\}}] dm \ = \ \int f d\nu_0^\alpha.
\end{eqnarray*}
If \(\int_E f(x) \mathbb{E}_x[e^{-\beta \zeta}1_{\{\partial\}}(X_{\zeta -})1_{\{s\geq \zeta\}}] dm = \infty\) for any \(s>0\), then we have
\[\int f d\nu_0^\beta \geq \int f d\nu_0^\alpha \geq \varlimsup_{s\searrow 0} \frac{e^{(\beta -\alpha) s}}{s} \int_E f(x) \mathbb{E}_x[e^{-\beta \zeta}1_{\{\partial\}}(X_{\zeta -})1_{\{s\geq \zeta\}}] dm = \infty.\]

\begin{rem}\label{remnu}
(1) The energy functional \(\nu_0\) is not a measure, but a functional. However, by \cite[Lemma 4.1]{Oo26}, the monotonicity and subadditivity hold for \(\nu_0\). By \cite[Proposition 4.2.3]{CF12}, it holds that \(\int |f|^2 d\nu_0=0\) for \(f\in \mathcal{F}\). For \(f\in C_0(E)\) and its compact support \(K\), by the regularity of a Dirichlet form, we can take a non-negative function \(g\in \mathcal{F} \cap C_0\) satisfying \(1_K \leq g\), so it holds that \(|\int f d\nu_0| \leq \int |f| d\nu_0 \leq \|f\|_{\infty} \int g^2 d\nu_0 = 0\). Hence the energy functional \(\nu_0\) can be regarded as a measure-like object supported at the cemetery point. 
Indeed, \(\nu_0\) plays the same role as the killing measure \(\kappa\) in the proofs, thus \(\kappa\) 
and \(\nu_0\) serve as parallel counterparts. See also \cite[Section 4]{Oo26} for details. \\
(2) We can replace \(\varlimsup_{s\searrow 0}\) with \(\lim_{s\searrow 0}\) in the above 
definition of \(\nu_0\) for an \(\alpha\)-excessive function \(f\), that is, a non-negative universally measurable function \(f\) satisfying \(e^{-\alpha t}p_t f \nearrow f\) as \(t \searrow 0\), because 
\(\frac{1}{s} \int_E f(x)(\varphi_{\alpha}(x)-e^{-\alpha s}p_s\varphi_{\alpha}(x)) dm(x)\) is increasing in 
\(s\searrow 0\).  This functional \(\int f d\nu_0\) was originally called an energy functional of \(f\) and 
\(\varphi_{\alpha}\) in \cite[Section 5.4]{CF12} and \cite[Section 3]{Ge90}.
\end{rem}

Similarly to Lemma \ref{Revuzkappa}, we have the following lemma.  We set \(\mathbb{E}_{\nu_0}[\cdot]:= \int_E \mathbb{E}_x[\cdot] d\nu_0(x).\)

\begin{lem}\label{Revuznu}
For \(\mu \in \mathcal{S}\), its corresponding PCAF \(\sfA\), any \(\alpha \geq 0\) and \(f \in \mathcal{B}_+(E)\), it holds that
\begin{equation}
\mathbb{E}_{\nu_0} \left[ \int_0^\infty e^{-\alpha s} f(X_s) \, d\sfA_s \right] = \int_E \varphi_{\alpha}f\, d\mu, \label{Revuz-nu-1}
\end{equation}
and, for \(t\geq 0\) such that there exists $\varepsilon > 0$ satisfying $\int_E \mathbb{P}_x(X_{\zeta-}=\partial, 
t<\zeta\le t + \vareps)  f(x) d\mu(x) < \infty$,
\begin{equation}
\mathbb{E}_{\nu_0} \left[ \int_0^t f(X_s) \, d\sfA_s \right] = \int_E (\varphi-p_t \varphi)f\, d\mu \label{Revuz-nu-2}.
\end{equation}
\end{lem}

\begin{proof}
First, we assume that \(\alpha >0\), \(\mu \in \mathcal{S}_0\) and \(f\) is a bounded positive Borel measurable function. Since \(\mathbb{E}_x \left[ \int_0^\infty e^{-\alpha s} f(X_s) \, d\sfA_s \right] \) is an \(\alpha\)-excessive quasi-continuous version of \(U_\alpha(f \mu)\), by Remark \ref{remnu}, we have
\[\mathbb{E}_{\nu_0} \left[ \int_0^\infty e^{-\alpha s} f(X_s) \, d\sfA_s \right]=
\lim_{s\searrow 0}\frac{1}{s} \int_E U_\alpha(f \mu)(x)(\varphi_{\alpha}(x)-e^{-\alpha s}p_s\varphi_{\alpha}(x)) dm(x)
 \]
and by applying \cite[Theorem 5.4.3.(iv)]{CF12} to the semigroup \(\{e^{-\alpha t}p_t\}_{t\geq 0}\), we have
\begin{eqnarray*}
\mathbb{E}_{\nu_0} \left[ \int_0^\infty e^{-\alpha s} f(X_s) \, d\sfA_s \right]= \int \varphi_{\alpha} f d\mu.
\end{eqnarray*}
Hence \eqref{Revuz-nu-1} holds for \(\mu \in \mathcal{S}_0\) and a bounded positive Borel measurable function \(f\). For any \(f \in \mathcal{B}_+(E)\), set \(f_n:=f\wedge n\). Then \(\{U_{\alpha}(f_n\mu)\}_n\) is an increasing sequence of \(\alpha\)-excessive functions, and 
\begin{equation}\label{eq:Revuznu-A}
\frac{1}{s}\int_E U_{\alpha}(f_n\mu)(x) (\varphi_{\alpha}(x)-e^{-\alpha s}p_s \varphi_{\alpha}(x))dm(x)
\end{equation}
is increasing in both \(s\) and \(n\) as \(s \searrow 0\) and \(n\nearrow \infty\). Hence, by taking the two limits in \eqref{eq:Revuznu-A} with respect to \(s\) and \(n\), which can be interchanged, \eqref{Revuz-nu-1} holds for \(\mu \in \mathcal{S}_0\) and \(f \in \mathcal{B}_+(E)\).

For any \(\mu \in \mathcal{S}\), we take a nest \(\{F_\ell\}\) such that \({1}_{F_\ell} \mu \in \mathcal{S}_0\) and \(\mu(\cap _\ell F_\ell^c)=0.\) By the same way as above, \eqref{Revuz-nu-1} holds for \(\mu \in \mathcal{S}\) and \(f \in \mathcal{B}_+(E)\).

For the case of \(\alpha =0\), we can prove \eqref{Revuz-nu-1} similarly to the above discussion using \(\mathcal{S}_0^{(0)}\) in place of \(\mathcal{S}_0\) due to \cite[Exercise 2.2.4]{FOT11} (see also the proof of \cite[Theorem 2.4.2.(ii)]{FOT11}).

Next we prove \eqref{Revuz-nu-2}. By the Revuz correspondence \eqref{eq-Revuz2}, we have
{
\begin{align}
\nonumber &\mathbb{E}_{\nu_0} \left[ \int_0^t f(X_u) \, d\sfA_u \right] = \varlimsup_{s \searrow 0} \frac{1}{s} \int_E \mathbb{E}_x \left[ \int_0^t f(X_u) \, d\sfA_u \right] (\varphi(x)-p_s\varphi(x)) dm(x)\\
&= \varlimsup_{s \searrow 0} \frac{1}{s} \int\int_0^t  p_u(\varphi - p_s \varphi)\, f\, du \,d\mu\ =\  \varlimsup_{s \searrow 0} \frac{1}{s} \int \int_0^t \mathbb{P}_x(X_{\zeta -}=\partial, u<\zeta\leq t+u)\, f(x)\, du \,d\mu\\
\nonumber &=\varlimsup_{s \searrow 0} \frac{1}{s} \int \mathbb{E}_x\left[\int_0^t 1_{\partial}(X_{\zeta -})1_{\{u<\zeta \leq s+u\}}\right]\,du\, f(x)\,d\mu\\
\nonumber &= \varlimsup_{s \searrow 0} \frac{1}{s} \int \mathbb{E}_x\left[1_{\partial}(X_{\zeta -}) (t\wedge \zeta - (\zeta-s)\vee0)\right]f(x)\,d\mu\\
&=\varlimsup_{s \searrow 0} \frac{1}{s} \int \mathbb{E}_x\left[1_{\partial}(X_{\zeta -}) (s\wedge \zeta - (\zeta-t)\vee0)\right]f(x)\,d\mu\ =\ \varlimsup_{s \searrow 0} \frac{1}{s} \int_0^s \int p_u(\varphi - p_t \varphi)\, f\, d\mu \,du \label{eq:Revuznu-C_A}\\
&\leq \varlimsup_{s \searrow 0} \int p_s(\varphi - p_t \varphi)\, f\, d\mu. \label{eq:Revuznu-C}
\end{align}
By \eqref{eq:Revuznu-C_A}, changing the variable, and Fatou's lemma, we have
\begin{align}
& \nonumber \mathbb{E}_{\nu_0} \left[ \int_0^t f(X_u) \, d\sfA_u \right] = \varlimsup_{s \searrow 0} \frac{1}{s} \int_0^s \int p_u(\varphi - p_t \varphi)\, f\, d\mu \,du\ \geq\  \varliminf_{s \searrow 0} \frac{1}{s} \int_0^s \int p_u(\varphi - p_t \varphi)\, f\, d\mu \,du\\
&= \varliminf_{s \searrow 0} \int_0^1 \int p_{us}(\varphi - p_t \varphi)\, f\, d\mu \,du\ \geq \ \int_0^1\int \varliminf_{s \searrow 0} p_{us}(\varphi - p_t \varphi)\, f\, d\mu \,du\\
& = \int_0^1\int \varliminf_{s \searrow 0} \mathbb{P}_x(X_{\zeta -}=\partial, us<\zeta\leq t+us)\, f\, d\mu \,du\ \geq \ \int_0^1\int \mathbb{P}_x(X_{\zeta -}=\partial, \zeta\leq t)\, f\, d\mu \,du\\
&= \int \mathbb{P}_x(X_{\zeta -}=\partial, \zeta\leq t)\, f\, d\mu \ = \ \int (\varphi - p_t \varphi)\, f\, d\mu.\label{eq:Revuznu-D}
\end{align}
}
\end{proof}

\begin{rem}
We remark that \(f\in L^1(\mu)\) satisfies the assumption in Lemma \ref{Revuznu}. Without the existence of \(\varepsilon >0\) satisfying \(\int_E \mathbb{P}_x(X_{\zeta -} = \partial, t< \zeta \leq t + \varepsilon)f(x)d\mu(x) < \infty\), the equality \eqref{Revuz-nu-2} does not hold in general. Let \(X\) be an absorbing Brownian motion on \((0,\infty)\), \(d\mu := dx\) and \(f(x) := e^{x^2/(2t)}x^{-1}(\log{x})^{-2}1_{\{x\geq e\}}\). By \cite[Theorem 2.35]{MP10}, we have
\[\mathbb{P}_x(\zeta \leq t) = 2\int_{\frac{x}{\sqrt{t}}}^\infty \frac{e^{-\frac{y^2}{2}}}{\sqrt{2\pi}}dy \leq \frac{2\sqrt{t}}{\sqrt{2\pi}x} e^{-\frac{x^2}{2t}},\] and so 
\[\int_E (\varphi-p_t \varphi)f\, d\mu  \leq \int_e^\infty \frac{2\sqrt{t}}{\sqrt{2\pi}x^2(\log{x})^2} dx = \frac{2\sqrt{t}}{e\sqrt{2\pi}} <\infty.\]
For \(s<(e^{-1}+t^{-1/2})^{-2}\), \(\frac{x}{\sqrt{s}} > \frac{x}{\sqrt{s+t}}+ 1\) holds for any \(x > e \), so, by \eqref{eq:Revuznu-C}, we have
\begin{align*}
&\mathbb{E}_{\nu_0} \left[ \int_0^t f(X_u) \, d\sfA_u \right] =  \varlimsup_{s \searrow 0} \int_e^\infty \int_{\frac{x}{\sqrt{s+t}}}^{\frac{x}{\sqrt{s}}}\frac{e^{-\frac{y^2}{2}}}{\sqrt{2\pi}}dy \frac{e^{\frac{x^2}{2t}}}{x(\log{x})^2} dx \\
&\geq  \varlimsup_{s \searrow 0} \int_e^\infty \int_{\frac{x}{\sqrt{s+t}}}^{\frac{x}{\sqrt{s+t}}+1}\frac{e^{-\frac{1}{2}(\frac{|x|}{\sqrt{s+t}}+1)^2}dy}{\sqrt{2\pi}} \frac{e^{\frac{x^2}{2t}}dx}{x(\log{x})^2}  \geq  \varlimsup_{s \searrow 0} \int_{e\vee \frac{1}{\sqrt{t+s}}\frac{4t(s+t)}{s}}^\infty  \frac{e^{\frac{sx^2}{2t(s+t)}-\frac{x}{\sqrt{s+t}}}e^{-\frac{1}{2}}}{\sqrt{2\pi}x(\log{x})^2}   dx \\
&\geq  \varlimsup_{s \searrow 0} \int_{e\vee \frac{1}{\sqrt{t+s}}\frac{4t(s+t)}{s}}^\infty  \frac{e^{\frac{sx^2}{4t(s+t)}}e^{-\frac{1}{2}}}{\sqrt{2\pi}x(\log{x})^2} dx  \ =\  \infty.
\end{align*}
\end{rem}

In \cite{Oo26}, the Revuz map from $({\cal S}_0, \rho_0)$ to $(\bfA_c^+, L^2(\mathbb{P}_{m+\kappa + \nu_0}))$ with the local uniform topology
is proved to be  a homeomorphism. This is accomplished by considering the 
relationship between PCAFs and smooth measures in the $L^2$ sense (see Corollary \ref{Oo26energy} below).
By considering this result with Theorem \ref{thm-Revuz}, Lemma \ref{Revuzkappa} and Lemma \ref{Revuznu},
we can establish a characterization of the Revuz correspondence in the $L^1$ sense as follows.

\begin{thm}\label{Revuz-trinity}
For \(\mu \in \mathcal{S}\), its corresponding  PCAF \(\sfA\), any \(\alpha \geq 0\) and \(f \in \mathcal{B}_+(E)\), it holds that
\begin{equation}
\mathbb{E}_{\alpha m+\kappa +\nu_0} \left[ \int_0^\infty e^{-\alpha s} f(X_s) \, d\sfA_s \right] = \int_E f\,  d\mu
\label{Revuz-trinity1}
\end{equation}
and, for \(t\geq 0\) such that there exists  \(\varepsilon >0\) satisfying \(\int_E \mathbb{P}_x(X_{\zeta -} = \partial, \zeta \leq t + \varepsilon)f(x)d\mu(x) < \infty\),
\begin{equation}
\mathbb{E}_{m} \left[ \int_0^t  f(X_s) \, d\sfA_s \right]+ \mathbb{E}_{\kappa+\nu_0} \left[\int_0^t (t-s) f(X_s) \, d\sfA_s \right] = t \int_E  f \, d\mu \label{Revuz-trinity2}
\end{equation}
\end{thm}

\begin{proof}
\eqref{Revuz-trinity1} follows immediately from \eqref{eq-Revuz} for \(h=1\),  \eqref{Revuz-kappa1} and \eqref{Revuz-nu-1}. 

By \eqref{Revuz-kappa2}, we have
\[\int_0^t \int_E ( 1_E-\varphi -p_u( 1_E-\varphi))\, f\, d\mu du = \int_0^t \mathbb{E}_{\kappa} \left [ \int_0^u  f(X_s)\, d{\sf A}_s \right]du = \mathbb{E}_{\kappa} \left[ \int_0^t (t-s)  f(X_s) \, d{\sf A}_s \right].\]

For \(t\geq 0\) such that there exists  \(\varepsilon >0\) satisfying \(\int_E \mathbb{P}_x(X_{\zeta -} = \partial, \zeta \leq t + \varepsilon)f(x)d\mu(x) < \infty\), \(u \leq t\) and \(s<\varepsilon\), as in \eqref{eq:Revuznu-C}, we have
\begin{align*}
\frac{1}{s}\int_E \mathbb{E}_{x} \left [ \int_0^u  f(X_v)\, d\sfA_v \right] (\varphi-p_s\varphi )\, dm &= \frac{1}{s} \int_0^s \int p_v(\varphi - p_u \varphi)\, f\, d\mu\,dv\\
&= \frac{1}{s} \int_0^s \int_E \mathbb{P}_x(X_{\zeta -} = \partial, v< \zeta \leq u + v)f(x)d\mu(x) \,dv\\
&\leq  \int_E \mathbb{P}_x(X_{\zeta -} = \partial, 0< \zeta \leq t + \varepsilon)f(x)d\mu(x),
\end{align*}
and so, by the dominated convergence theorem and the proof of Lemma \ref{Revuznu},
\[\int_0^t \mathbb{E}_{\nu_0} \left [ \int_0^u  f(X_v)\, d\sfA_v \right]du = \lim_{s\searrow 0} \int_0^t  \frac{1}{s}\int_E \mathbb{E}_{x} \left [ \int_0^u  f(X_v)\, d\sfA_v \right] (\varphi-p_s\varphi ) \, dm\, du. \]

Combining this with \eqref{Revuz-nu-2} and the monotone convergence theorem, we have
\begin{eqnarray*}
\int_0^t \int_E (\varphi -p_u\varphi)\,f\,  d\mu du &=& \int_0^t \mathbb{E}_{\nu_0} \left [ \int_0^u  f(X_v)\, d\sfA_v \right]du \\
&=& \lim_{s\searrow 0} \frac{1}{s}\int_E \mathbb{E}_{x} \left [ \int_0^t \int_0^u  f(X_v)\, d\sfA_vdu \right] (\varphi-p_s\varphi )\, dm\\
&=& \lim_{s\searrow 0} \frac{1}{s}\int_E \mathbb{E}_{x} \left [ \int_0^t (t-v) f(X_v)\, d\sfA_v \right] (\varphi-p_s\varphi ) \, dm\\
&=& \mathbb{E}_{\nu_0} \left [ \int_0^t (t-v) f(X_v)\, d\sfA_v \right],
\end{eqnarray*}
so \eqref{Revuz-trinity2} holds.
\end{proof}

By using Theorem \ref{Revuz-trinity}, we can obtain the following corollary by another way of \cite[Proposition 1.2]{Oo26} of independent interest.
\begin{cor}[{\cite[Proposition 1.2]{Oo26}}] \label{Oo26energy}
For \(\mu, \nu \in \mathcal{S}_0\) and their corresponding  PCAFs \(\sf A, \sf B\), it holds that
\begin{equation}
\mathbb{E}_{\alpha m+\frac{\kappa}{2}+\frac{\nu_0}{2}}\left[\widetilde{{\sf A}_{\infty}}\widetilde{{\sf B}_{\infty}}\right] =\mathcal{E}_{\alpha}(U_{\alpha}\mu, U_{\alpha}\nu),
\end{equation}
where we set \(\widetilde{{\sf A}_t}:=\int_0^t e^{-\alpha s}d\sf A_s\) and \(\widetilde{{\sf B}_t}:=\int_0^t e^{-\alpha s}d\sf B_s\) for \(\alpha>0\). 
\end{cor}
\begin{proof}
 By the polarization identity, it suffices to consider the case \(\sf A =\sf B\) and, by \cite[Theorem 2.2.4]{FOT11}, we may assume that \(U_{\alpha}\mu\) is bounded. By the resolvent equation and \cite[Exercise 4.1.7]{CF12}, we have
 \begin{align*}
  U_{\alpha}(U_{\alpha}\mu\cdot \mu)(x) &=
 U_{2\alpha}(U_{\alpha}\mu\cdot \mu)(x) +
 \alpha R_{\alpha}U_{2\alpha}(U_{\alpha}\mu\cdot \mu)(x) \\
 &= \frac{1}{2}\mathbb{E}_x \left[(\widetilde{{\sf A}_{\infty}})^2\right]
 +  \frac{\alpha}{2}R_{\alpha}\mathbb{E}_{\cdot}
 \left[(\widetilde{{\sf A}_{\infty}})^2\right](x)
 \end{align*}
and by Theorem \ref{Revuz-trinity}, we have
\begin{eqnarray*}
\lefteqn{ \mathcal{E}_{\alpha}(U_{\alpha}\mu, U_{\alpha}\mu) = \int_E U_{\alpha}\mu \,d\mu \ =\  \mathbb{E}_{\alpha m+\kappa +\nu_0} \left[ \int_0^\infty e^{-\alpha s} U_{\alpha}\mu (X_s) \, d\sfA_s \right]}\\
&=& \int_E U_{\alpha}(U_{\alpha}\mu\cdot \mu)\, d(\alpha m+\kappa +\nu_0)\\
&=& \mathbb{E}_{\alpha m+\frac{\kappa}{2}+\frac{\nu_0}{2}} \left[(\widetilde{{\sf A}_{\infty}})^2\right] -\frac{\alpha}{2} \mathbb{E}_m\left[ (\widetilde{{\sf A}_{\infty}})^2 \right] + \frac{\alpha}{2} \int_E R_{\alpha}\mathbb{E}_{\cdot} \left[ (\widetilde{{\sf A}_{\infty}})^2 \right]\, d(\alpha m +\kappa +\nu_0)\\
&=& \mathbb{E}_{\alpha m+\frac{\kappa}{2}+\frac{\nu_0}{2}} \left[(\widetilde{{\sf A}_{\infty}})^2\right].
\end{eqnarray*}
Here, in the last equality, for $\dis f(x)=\mathbb{E}_x \left[(\widetilde{{\sf A}_{\infty}})^2\right],$ we used
$\dis \int \alpha R_{\alpha}f dm =\int \alpha R_{\alpha}1 \cdot f dm$, 
$$\dis\int R_{\alpha}f d\kappa =\int (1-\alpha R_{\alpha}1-\varphi_{\alpha}) f dm
\quad {\rm  and} \quad \dis \int R_{\alpha}f d\nu_0 =\int \varphi_{\alpha} f dm.
$$
\end{proof}

Using Theorem \ref{Revuz-trinity}, we obtain the following sufficient condition for vague convergence of smooth measures. 
\begin{thm}\label{TVconv}
For Radon measures \(\mu_n, \mu \in \mathcal{S}\) and their corresponding PCAFs \({\sf A}^n, {\sf A}\), we assume that \({\sf A}^n\) converges to \({\sf A}\) in \(L^1(\mathbb{P}_{m+\kappa+\nu_0})\) with the local total variation distance, that is, for any \(T>0\),
\[\lim_{n\to \infty} \mathbb{E}_{m+\kappa+\nu_0}\left[TV_{[0,T]}({\sf A}^n - {\sf A}) \right]=0,\]
where \(TV_{[0,T]}(\cdot)\) denotes the total variation on \([0,T]\). Then \(\mu_n\) converges vaguely to \(\mu \).
\end{thm}

\begin{proof}
Assume that \(\lim_{n\to \infty} \mathbb{E}_{m+\kappa+\nu_0}\left[TV_{[0,T]}({\sf A}^n - {\sf A}) \right]=0\). Then there exists \(C_T>0\) such that \(\sup_n \mathbb{E}_{m+\kappa+\nu_0}\left[TV_{[0,T]}({\sf A}^n - {\sf A})\right] \leq C_T\). For any \(t\) and \(f\in \mathcal{F} \cap C_0(E)\), it holds that \(\int_E \mathbb{P}_x(X_{\zeta -} = \partial, t< \zeta \leq t + \varepsilon)f(x)d\mu_n(x) \leq \int fd\mu_n < \infty\), and so, by Theorem \ref{Revuz-trinity},  we have
\begin{eqnarray*}
t \left| \int_E  f \, d\mu_n -\int_E  f \, d\mu \right| &\leq & \mathbb{E}_{m} \left| \int_0^t  f(X_s) \, d({\sf A}_s^n-{\sf A}_s) \right|+ \mathbb{E}_{\kappa +\nu_0} \left|\int_0^t (t-s) f(X_s) \,d({\sf A}_s^n-{\sf A}_s) \right|\\
&\leq & \|f\|_{\infty}  \mathbb{E}_{m} \left[ TV_{[0,t]}({\sf A}^n-{\sf A}) \right]+ t \|f\|_{\infty}\mathbb{E}_{\kappa+\nu_0} \left[TV_{[0,t]}({\sf A}^n-{\sf A}) \right]\\
&\leq & (1+t)\|f\|_{\infty}  \mathbb{E}_{m+\kappa+\nu_0} \left[ TV_{[0,t]}({\sf A}^n-{\sf A}) \right].
\end{eqnarray*}
As \(n\) tends to infinity, \(\int_E  f \, d\mu_n\) converges to \(\int_E  f \, d\mu\) for \(f\in \mathcal{F} \cap C_0\).
For any \(f\in C_0(E)\), we take \(f_{\varepsilon} \in \mathcal{F} \cap C_0(E)\) satisfying \(\|f-f_{\varepsilon}\|_{\infty}\leq \varepsilon\) and supp\((f_{\varepsilon}) \subset \) supp\((f)=:K\). We take a non-negative function \(\varphi \in \mathcal{F} \cap C_0(E)\) satisfying \(\varphi =1\) on \(K\), then we have
\begin{eqnarray*}
\left|\int f_{\varepsilon}d\mu_n -\int fd\mu_n \right| &\leq & \varepsilon \mu_n(K)\ \leq \ \varepsilon \int \varphi d\mu + \varepsilon \left|\int \varphi d\mu_n-\int \varphi d\mu \right|\\
&\leq &\varepsilon \int \varphi d\mu + \varepsilon (1+\frac{1}{t})\|\varphi\|_{\infty} \mathbb{E}_{m+\kappa +\nu_0} \left[ TV_{[0,t]}({\sfA^n}-{\sfA}) \right]\ \leq \  \varepsilon C
\end{eqnarray*}
where \(C:=\left( \int\varphi d\mu + C_t\left(1+ \frac{1}{t}\right) \right)\) is independent of \(n.\) Combining this with the convergence of \(\int f_{\varepsilon} d\mu_n\), the integral \(\int f d\mu_n\) converges to \(\int f d\mu\). 
\end{proof}

\begin{thm}
Suppose \(\mathbb{P}_x(X_{\zeta -} \in E, \zeta <\infty) =  \mathbb{P}_x(\zeta <\infty)\) for q.e. \(x\) and \(m(E)+\kappa(E)<\infty \).  For Radon measures \(\mu_n, \mu \in \mathcal{S}\) and their corresponding PCAFs \(\sf A^n, A\), if \(\sf A^n\) converges to \(\sf A\) in \(L^1(\mathbb{P}_{m+\kappa})\) with the local uniform topology, then \(\mu_n\) converges vaguely to \(\mu\).
\end{thm}
\begin{proof}
We remark that \(\nu_0 =0\) follows from the assumption \(\mathbb{P}_x(X_{\zeta -} \in E, \zeta <\infty) =  \mathbb{P}_x(\zeta <\infty)\).

Set \(D(\mathbb{R}) := \{ f: [0,\infty) \to \mathbb{R} \mid  f \text{ is c\`adl\`ag} \}\) equipped with the local uniform topology. By \cite[Section VI.1.1a]{JS03}, this space is a complete metric space but not separable. Set \(C_{\uparrow}(\mathbb{R}_+):=\{f : [0,\infty) \to \mathbb{R} \mid f \text{\ is\ non-decreasing\ and\ continuous\ with\ }f(0)=0. \}\) and define $\Phi : D(\mathbb{R})\times C_{\uparrow}(\mathbb{R}_+) \to \mathbb{R}$ by
\[
\Phi(f, a) = \int_0^t f(s) \, da(s)
\]
for fixed \(t\). If \((f_n, a_n)\) converges to \((f, a)\) with the local uniform topologies, then we have
\begin{eqnarray*}
|\Phi (f_n, a_n)- \Phi(f,a)| & \leq & \int_0^t |f_n(s)-f(s)|da_n(s) +\left|\int_0^t f(s) (da_n(s)-da(s)) \right|\\
&\leq & \sup_{s\leq t} |f_n-f| \cdot t\cdot \sup_{s\leq t}|a_n(s)| +\left|\int_0^t f(s) (da_n(s)-da(s)) \right|
\end{eqnarray*}
and the second term converges to \(0\) by \cite[Lemma 4.8]{Oo25}, so \(\Phi\) is a continuous function. Since \(\sf A^n\) converges to \(\sf A\) in \(L^1(\mathbb{P}_{m+\kappa})\) and \(\mathbb{P}_{m+\kappa}\) is a finite measure, \(\sup_{0\leq s\leq t}|{\sf A}_s^n - {\sf A}_s|\) converges to \(0\) in measure \(\mathbb{P}_{m+\kappa}\) and we have that \(\{{\sf A}_t^n\}_n\) is uniformly integrable under \(\mathbb{P}_{m+\kappa}\). Hence, by the continuous mapping theorem for probability measures \cite[Theorem 2.3.4]{D10}, \(\Phi(f(X_{\cdot}), \sf A^n)\) converges to \(\Phi(f(X_{\cdot}), \sf A) \) in measure \(\mathbb{P}_{m+\kappa}\) for any \(f\in C_0(E).\) By the uniform integrability of \(\{{\sf A}_t^n\}_n\), we have
\[\sup_n \mathbb{E}_{m+\kappa}\left|\int_0^tf(X_s)d\sfA_s^n \right| \leq \sup_n \|f\|_{\infty} \mathbb{E}_{m+\kappa}[{\sf A}_t^n] <\infty. \]
For any \(\varepsilon>0,\) there exists \(\delta>0\) such that, if \(\mathbb{P}_{m+\kappa}(\Lambda) \leq \delta\) then \( \sup_n \mathbb{E}_{m+\kappa}[\sfA_t^n 1_{\Lambda}]\leq \varepsilon .\) Hence we have
\[\sup_n \mathbb{E}_{m+\kappa}\left|\int_0^tf(X_s)d\sfA_s^n 1_{\Lambda} \right| \leq \sup_n \|f\|_{\infty} \mathbb{E}_{m+\kappa}[\sfA_t^n 1_{\Lambda}]\leq \varepsilon \|f\|_{\infty} ,\]
 and so \(\{ \int_0^tf(X_s)d\sfA_s^n \}_n\) is uniformly integrable under \(\mathbb{P}_{m+\kappa}\) and \(\int_0^tf(X_s)d\sfA_s^n\) converges to \(\int_0^tf(X_s)d{\sf A}_s\)
 in \(L^1(\mathbb{P}_{m+\kappa})\).
Since it holds that, for \(T\geq 0\),
\begin{eqnarray*}
\lefteqn{ \mathbb{E}_{\kappa} \left[\sup_{0\leq t\leq T} \left|\int_0^t (t-s)  f(X_s) \, d{\sf A}_s^n-  \int_0^t (t-s)  f(X_s) \, d{\sf A}_s \right|\right]}\\
&=& \mathbb{E}_{\kappa}  \left [\sup_{0\leq t\leq T} \left|\int_0^t \int_0^u  f(X_s)\, d{\sf A}_s^ndu-\int_0^t \int_0^u  f(X_s)\, d{\sf A}_s du \right|\right]\\
&\leq & \mathbb{E}_{\kappa}  \left [T\sup_{0\leq u\leq T}\left| \int_0^u  f(X_s)\, d{\sf A}_s^n- \int_0^u  f(X_s)\, d{\sf A}_s \right|\right],
\end{eqnarray*} 
by Theorem \ref{Revuz-trinity}, \(\mu_n\) converges vaguely to \(\mu\).
\end{proof}

\hspace{-7mm}\textbf{Acknowledgements}\\
The authors would like to thank Prof. Masayoshi Takeda for useful discussions and continuous encouragement.  
His helpful comments on the proof of Lemma \ref{Revuzkappa} in an earlier version of this paper were particularly valuable. This work was supported by JSPS KAKENHI Grant Numbers  25K17270 (T.O.) and 25K07056 (T.U.).

\end{document}